\documentclass[12pt]{amsart}
\usepackage{float}
\usepackage{amsmath}
\usepackage{amssymb}
\usepackage{amsthm}
\usepackage{amsfonts}
\usepackage[all]{xy}
\usepackage[margin=3cm]{geometry}
\usepackage{mathtools}
\usepackage{qtree}
\usepackage{tikz}
\usepackage[colorlinks=true]{hyperref}
\usepackage{mathtools}
\usepackage{enumerate}
\usepackage{pifont}
\usepackage{wasysym}
\usepackage{stmaryrd}
\usepackage{theoremref}
\newcommand{\CC}{\mathbb{C}}

\newcommand{\PP}{\mathbb{P}}

\newcommand{\iso}{\cong}

\newcommand{\wbuild}{\mathcal{G}_{\mathcal{A}}}
\newcommand{\wbuildb}{\mathcal{G}_{\mathcal{B}}}
\newcommand{\wbuildc}{\mathcal{G}_{\mathcal{C}}}
\newcommand{\wmod}{X_\mathcal{A}[n]}

\newcommand{\wmodb}{X_\mathcal{B}[n]}
\newcommand{\wmodc}{X_\mathcal{C}[n]}

\newcommand{\wbuilds}{\mathcal{G}_{\mathcal{A'}}}
\newcommand{\wmods}{X_\mathcal{A'}[r]}

\newcommand{\wbuildrel}{\mathcal{G}_{\mathcal{A}}^{rel}}
\newcommand{\wbuildrelu}{\mathcal{G}_{\mathcal{A}}^{rel,+}}
\newcommand{\wbuildrelU}{\mathcal{G}_{\mathcal{A},U}^{rel}}
\newcommand{\wmodrel}{X_Y^{\mathcal{A}}[n]}
\newcommand{\ch}{\mathcal{C}_{N,S}}
\newcommand{\cc}{\mathcal{C}_{N,S}'}

\newtheorem{myex}{Example}
\newtheorem{mydef}{Definition}
\newtheorem{thm}{Theorem}
\newtheorem{lem}{Lemma}

\newtheorem{prop}{Proposition}

\title[ Weighted compactifications of configuration spaces]{ Weighted compactifications of configuration spaces and relative stable degenerations}
\begin{document}
\author{Evangelos Routis}
 \thanks{webpage: \url{http://www.math.brown.edu/~r0utis/}}
\address{Brown University, Department of Mathematics,
Box 1917,
151 Thayer Street
Providence, RI 02912
}
\email{r0utis@math.brown.edu}

\maketitle

\begin{abstract}We study a compactification of the configuration space of $n$ distinct labeled points on an arbitrary nonsingular variety. Our construction provides a generalization of the original Fulton-MacPherson compactification that is parallel to the generalization of the moduli space of $n$-pointed stable curves carried out by Hassett.
\end{abstract}

\begin{section}{Introduction}

As configurations of points on a projective curve come up naturally in many contexts,  moduli spaces $\operatorname{M}_{g,n}$, of smooth $n$ pointed  curves of genus $g$ arise naturally as a means by which to study them.  While not themselves compact, these varieties do admit  a number of compactifications.  For example,  by the moduli space $\overline{\operatorname{M}}_{g,n}$ of $n$-pointed stable curves  due to Deligne, Mumfurd and Knudsen (\cite{DM}, \cite{KnM},\cite{Kn2}, \cite{Kn3}), and also by 
moduli spaces $\overline{M}_{g,\mathcal{A}}$, due to Hassett \cite{Hassett}, parametrizing weighted configurations of  $n$-pointed stable curves of genus g.  

 Hassett's weighted spaces receive birational morphisms from 
$\overline{\operatorname{M}}_{g,n}$, and these maps give useful tools for understanding the spaces $\overline{\operatorname{M}}_{g,n}$. Recently these maps were used to settle a 15 year old question about  $\overline{\operatorname{M}}_{0,n}$, about whether or not it is a so-called Mori Dream Space (MDS). A special case of Hassett spaces are known to be toric varieties (the Losev Manin spaces $\overline{LM}_n$) , and as such are MDS. In recent celebrated work,  Castravet and Tevelev \cite{C-T} prove, that for $n=134$, the blowup $BL_{e}\overline{LM}_{134}$ of the Losev Manin Spaces at the identity of its torus, is not a Mori Dream Space, and so consequently as an MDS must map onto an MDS, $\overline{M}_{0,n}$ is not an MDS, for $n \ge 134$.

   Such weighted spaces arise in other contexts.   Alexeev and Guy (\cite{AG}) and, independently, Bayer and Manin (\cite{BayerManin}) introduced the theory of weighted stable maps, extending Kontsevich's theory of ordinary stable maps. Special cases of weighted stable maps also appear in the work of A. Musta\c{t}\u{a} and M-A. Musta\c{t}\u{a} \cite{MM}.

In this work we introduce and study a nonsingular variety $\wmod$, whose points correspond to stable weighted configurations of points on a smooth variety $X$.  We show that  $\wmod$ compactifies the configuration space $F(X,n)$, of $n$ distinct labeled points on $X$, and generalize  the Fulton 
MacPherson spaces $X[n]$. We describe the boundary, showing it is a normal crossings divisor (Theorem 2). Just as with the Hassett spaces, we show that the $\wmod$ receive birational morphisms from $X[n]$ consisting of a sequence of blowups (Theorem \ref{thm5}). Our construction generalizes to a relative context, which also gives a relative construction for the original spaces $X[n]$.  As a first application, we introduce the moduli space of \textit{n pointed relative} $\mathcal{A}$ \textit{stable} \textit{degenerations}(Section \ref{sec3}). As a second application, we provide a presentation of the Chow ring of $\wmod$(Theorem \ref{thm8}), which in turn can be used to compute the Chow ring of Hassett's space in genus 0\cite{R2}.
\tableofcontents
\subsection{Precise statement of results}

Let $X$ be a nonsingular variety over an algebraically closed field $k$  and $\mathcal{A}:=\{a_1,a_2, \dots,a_n\} $ be a set of rational numbers such that\\
\renewcommand{\labelitemi}{$\bullet$}
\begin{itemize}
\item $0<a_i\leq 1 , i=1,2,\dots,n$ 
\end {itemize}
\vspace{0.1in}
The aim of this paper is to introduce and study the so called \textit{weighted compactification} $\wmod$ of the configuration space $F_\mathcal{A}(X,n)$ of $n$ labeled points in $X$ carrying weights $a_i$, that is the parameter space of ordered $n$-tuples of points equipped with weights $a_i$ that satisfy the condition:\\
\renewcommand{\labelitemi}{$\bullet$}
\begin{itemize}
\item for any set of labels $S\subset N$ of coincident points we have $\sum \limits _{i\in S} a_i \leq 1\ $.
\end {itemize}
\vspace{0.1in}
The resulting space is a generalization of the Fulton Macpherson compactification (\cite{FM}) and also provides a compactification of the configuration space of $n$ \textit{distinct} labeled points. The main definitions and properties are found in sections \ref{sec2.2} and \ref{sec2.3}. Let

\begin{center}
$\mathcal{G}:=\{ \Delta_{12\dots n}; \Delta_{12\dots (n-1)},\dots ,\Delta_{23\dots n}; \dots;  \Delta_{12},\dots   ,\Delta_{(n-1)n} \} $
\end{center}

be the set of diagonals of $X^n$ listed in ascending dimension order.
\vspace{0.1in}

Also, let $N:=\{1,2,\dots,n\}$ and set
\begin{center}
$\mathcal{G}_{\mathcal{A}}:= \{ \Delta_S\subset X^n | S\subset N \,\text{and}\, \sum \limits _{i\in S} a_i > 1\}$
\end{center}
 and list its elements in ascending dimension order. \\
 
\textbf{Definition}: The \textit{weighted compactification} $\wmod$ of the configuration space $F_\mathcal{A}(X,n)$ is the closure of the image of the diagonal embedding\\
 
  \begin{center}
 
$ X^n\setminus \bigcup\limits_ {\Delta_I \in \mathcal{G}_{\mathcal{A}}} \Delta_I\xhookrightarrow\ X^n \times\displaystyle \prod \limits_ {\Delta_I\in\ \mathcal{G}_{\mathcal{A}}} Bl_{\Delta_I} X^n$

\end {center}
\vspace{0.1in}
\textbf{Main Theorem:}(Section \ref{sec2.2}, Theorem 2)
\begin{enumerate}
\item{$\wmod$ is a nonsingular variety. The boundary  $X_\mathcal{A}[n] \setminus (X^n \setminus \bigcup \limits_ {\Delta_I \in \mathcal{G}_{\mathcal{A}}} \Delta_I) $ is the union of $|\mathcal{G}_{\mathcal{A}}|$ divisors $D_I,$ where $I\subset N   , |I|\geq2$ and $ \sum \limits _{i\in I} a_i >1\ $.}

\item{ Any set of the boundary divisors intersects transversally. An intersection of divisors $D_{I_1}\cap D_{I_2}\cap\dots D_{I_k}$ is nonempty precisely when the sets are nested in the sense that any pair $\{I_i,I_j\}$ either has empty intersection or one set is contained in the other. }

 \item{$X_\mathcal{A}[n]$ is the iterated blowup of $X^n$ along dominant transforms of the diagonals of $\mathcal{G}_{\mathcal{A}}$} in the above order.
\end{enumerate}
Note: The \textit{dominant transform} of a variety under a blowup is either the strict transform, in case the variety is not contained in the center of the  blowup, or the inverse image of the variety via the blowup map if it is contained in the center.\\

We further introduce the \textit{universal family} $\wmod^+\rightarrow \wmod$, in section \ref{sec2.3}. This is a flat morphism between nonsingular varieties equipped with $n$ sections in its smooth locus (cf Theorem \ref{thm3}). Its fibers, called $n$-\textit{pointed} $\mathcal{A}$ \textit{stable degenerations}, are described in section \ref{sec2.2}. These look like the degenerations that one sees in the Fulton Macpherson compactification, with the exception that the sections now carry weights and are allowed to collide depending on the accumulation of their weights. We also provide a functorial description of the weighted compactification in section \ref{sec2.4} in analogy with the one in \cite{FM}. Forgetful and reduction morphisms analogous to the ones in \cite{Hassett}, \cite {BayerManin}, \cite{AG} are introduced in section \ref{sec2.5}. \\ 

In \cite{FM}, Fulton and Macpherson state that one should be able to extend their results by providing a compactification of the configuration space of $n$ distinct labeled points on the fibers of a smooth family $X\rightarrow Y$.  Pandharipande (\cite{Pa1}) shows that this is possible under the additional condition that the morphism $X\rightarrow Y$ is projective. We show that the Fulton-Macpherson construction can be lifted to the relative case without the projectivity assumption and more generally, that there exists a moduli space $X_Y^{\mathcal{A}}[n]$ that parametrizes $n$ pointed $\mathcal{A}$ stable degenerations on the fibers of $X\rightarrow Y$. We call the latter the \textit{relative weighted compactification} (Section \ref{sec3.1}). As expected, the properties described in Section 2 regarding the absolute case $X\rightarrow k$ carry over to the relative setting (see Section \ref{sec3.1} for a discussion).  As an application,  we construct the moduli space of $n$-pointed relative $\mathcal{A}$ stable degenerations and its universal family (section \ref{sec3.2}):  let $W$ be a smooth variety and $B$ a smooth curve with a distinguished point $b_0 \in B$. Also, let $\pi: W\rightarrow B$ be a flat morphism  such that the fibers $W_b$ are smooth for all $b \in B$ except  $b=b_0$, in which case $W_{b_0}$ is the union of two smooth varieties $X_1$ and $X_2$ intersecting transversally along a smooth divisor $D$. This implies that $N_{D/ X_1}\iso N_{D/X_2}^{\vee}$. We have the following :\\

 \textbf{Definition} (\cite[Definition 2.3.1] {ACFW}):  The \textit{expansion of length} $l\geq 0$ of $W_{b_0}$ is  the variety
     \begin{gather*}
     W(l):=X_1\bigsqcup\limits_{D=D_1^-} P_1\bigsqcup_{D_1^+=D_2^-}\cdots \bigsqcup_{D_{l-1}^+=D_l^-}P_l  \bigsqcup_{D_{l+}^+=D}X_2
               \end{gather*}
        where the $P_i$ are isomorphic to the $\PP^1$- bundle $\PP(N_{D/X_1}\oplus \mathcal{O}_D)=\PP( \mathcal{O}_D \oplus N_{D/X_2})$ and $D_j^-, D_j^+$ are the zero and infinity sections of $P_j$ respectively. The $P_j$ are called the exceptional components of the expansion.
     \bigskip
    
     In \cite{AF2}, the authors construct a compactification of the space of $n$ tuples of points on the fibers of $\pi$ such that the points are not allowed to land in the singular locus, by replacing the critical fibers is by expansions. We give the basic definitions and refer the reader to [ibid.] for more details. \\
     
\textbf{Definition}(\cite[{Definition 1.2.1}]{AF2}): A \textit{stable expanded configuration} $(\mathcal{W},\widetilde{\sigma_i})$ of degree $n$ on $\pi:W\rightarrow B$ consists of:
     \begin{enumerate}
     \item either a smooth fiber $W_b$ or an expansion of $W_{b_0}$, which we denote by $\mathcal{W}$
     \item an ordered collection of $n$ smooth points $\widetilde{\sigma_i}\in \mathcal{W}^{sm},\,i=1,\dots n$
     \end{enumerate}
     such that in case $\mathcal{W}$ is an expansion, the following stability condition holds:
     \begin{itemize}
     \item each exceptional component $P_j$ contains at least one of the $\widetilde{\sigma_i}$.
     \end{itemize}
     \bigskip

 Let $\mathcal{A}:=\{a_1,a_2, \dots,a_n\} $ be a set of rational numbers as above. \\
 
   \textbf{Definition}: An $n$-\textit{pointed relative} $\mathcal{A}$ \textit{stable degeneration} of $W\rightarrow B$ is pair $(\mathcal{W'},\sigma_i')$ such that 
      \begin{center}
    $W'=\mathcal{W}\bigsqcup\limits_{\mathcal{W}^{sm}\setminus\{\sigma_i\}}\widetilde{ \mathcal{W}}$
    \end{center}
    where $(\mathcal{W},\sigma_i)$ is a stable expanded configuration of degree $n$ and $(\widetilde{ \mathcal{W}},\sigma_i')$ is an $n$-pointed $\mathcal{A}$ stable degeneration of $(\mathcal{W}^{sm},\sigma_i)$ in the sense of section \ref{sec2.2}.    \\
    
The moduli space that parametrizes such objects is discussed in section \ref{sec3.2}. The motivation for introducing such degenerations is the construction of the moduli space of \textit{relative (un)ramified stable maps}, that is maps from curves to $n$ pointed relative stable degenerations (with all weights equal to 1). This is the relative version of the moduli space of (un)ramified stable maps \cite{KKO}, which in turn generalizes the moduli space of admissible covers (\cite {H-M}, \cite{ACV}), by allowing target degenerations to be the Fulton Macpherson degenerations. This project will be carried out elsewhere.\\

In section \ref{sec4}, we provide some further relevant examples from the literature and show how they relate to our construction. \\

One  the advantages of working with nonsingular spaces is that one often has enough tools to determine their intersection theory. The blowup presentation of $\wmod$ is an asset that makes things simpler. By making use of the explicit geometry of the strict transforms of the diagonals in $\wbuild$ in each blowup stage, as demonstrated in detail in Section \ref{sec5}, and by successively applying Keel's blowup formula (\cite{Keel}), we conclude the complete description of the Chow ring of $\wmod$. 
We assume that  $X$ has a \textit{cellular decomposition}, i.e. there exists a filtration\\
\begin{gather*}
X_n:=X\supset X_{n-1}\supset \dots X_0  \supset X_{-1}:=\emptyset 
\end{gather*}
such that each $X_i\setminus X_{i-1}$ is isomorphic to a disjoint union of affine spaces. The reason we make this assumption is that in this case there exists a natural isomorphism $A^{\bullet}(X^n)\iso (A^{\bullet}X)^{\otimes n}$. \\

Now, let $S\subset N$ such that $\Delta_S\in \wbuild$. The restriction morphism $A^{\bullet}X^n\rightarrow A^{\bullet}\Delta_S$ is surjective. Define \\
\begin{gather*}
J_S:=ker(A^{\bullet}X^n\rightarrow A^{\bullet}\Delta_S)
\end{gather*}
By our assumptions, $X^n$ has a K\"{u}nneth decomposition, so we have\\
\begin{gather*}
J_S=\left(p_a^*x-p_b^*x\right)_{x\in A^{\bullet}X,\, a,b \in S}
\end{gather*}
where $p_a:X^n\rightarrow X$ is the projection onto the $a$-th factor. Also, define the polynomial $c_{a,b}(t) \in A^{\bullet}(X^n)[t]$ by \\
\begin{gather*}
c_{a,b}(t):=\displaystyle \sum\limits_{i=1}^{m}(-1)^i p_a^*(c_{m-i}) t^i +[\Delta_{\{a,b\}}]
\end{gather*}
where $m=dim (X)$ and $c_{m-i}$ is the $(m-i)$th Chern class of the tangent bundle $T_X$ of $X$. We give the following:\\

\textbf{Definition:}
\begin{enumerate}
\item Let $S,T\subset N$ such that $|S|,|T|\geq2$. We say that $S$ and $T$ \textbf{overlap} if $S$ and $T$ have nonempty intersection and neither of them  is contained in the other. In this case we also say that $T$ is an \textbf{overlap} of $S$ and vice versa.
\item Let $S,T\subset N$ such that $|S|,|T|\geq2$. We say that $T$ is a \textbf{weak overlap} of $S$ (and vice versa) if the intersection of $S$ and $T$ is a singleton.
\end{enumerate}
\textbf{Remark}: In the above definition we do not require $\Delta_S$ or $\Delta_T$ to be in $\wbuild$.\\

We are now ready to state the other main result of this paper: let $A^{\bullet}(X^n)[D^S]$ be the polynomial ring over $A^{\bullet}(X^n)$ in the variables $D^S$, where $S\subset N$ is such that $\Delta_S\in \wbuild$. This polynomial ring maps to the Chow ring of $\wmod$ by sending each $D^S$ to the class of the divisor $D_S$ of $\wmod$.\\

\textbf{Main Theorem 2} (Section 5, Theorem \ref{thm8}):With the above assumptions and notation, the Chow ring of $\wmod$ is the quotient
\begin{gather*}
\displaystyle\dfrac{A^{\bullet}(X^n)[D^S]}{\mathcal{I}}
\end{gather*}
where $\mathcal{I}$ is the ideal generated by the following elements:
\begin{itemize}
\item $D^S\cdot D^T$, for any $S$ and $T$ that overlap and both $\Delta_S$ and $\Delta_T$ belong to $\wbuild$;\\
\item $J_S\cdot D^S$, where $\Delta_S\in \wbuild$;\\
\item $D^S\cdot \displaystyle \prod\limits_{k=1}^{|T|-1} c_{i_k,i_{k+1}}\left(\sum \limits_{I\supseteq S\cup T} D^I\right)$, where $T=\{i_1, i_2,\dots ,i_{|T|}\}$ is a weak overlap of $S$, $\Delta_S\in \wbuild$ and $\Delta_T$ is not necessarily in $\wbuild$;\\
\item$\displaystyle \prod\limits_{k=1}^{|S|-1} c_{i_k,i_{k+1}}\left(\sum \limits_{I\supseteq S} D^I\right)$, where $S=\{i_1,i_2,\dots ,i_{|S|}\}$ and $\Delta_S\in \wbuild$.
\end{itemize}
\vspace{0.1in}

 In \cite{R2} we use the results of this last section to compute the Chow ring of Hassett's space in genus 0 \cite{Hassett}. The advantage of this approach is that it avoids all GIT constructions of the moduli space of weighted stable curves as presented for example in \cite{KM} and \cite{GJM}. The construction of Kiem and Moon is given via a sequence of blowups starting from $(\PP^1)^n \sslash SL_2$  and exhibits $\overline{M}_{0,\mathcal A}$ as  a GIT quotient of the moduli space of parametrized rational stable curves (\cite{MM}) by $SL_2$ for special cases of $\mathcal A$. However, the blowup centers are not always smooth. In fact, when $n$ is even the above sequence involves Kirwan's blowup desingularization. \\

\textbf{Conventions}: Throughout this paper we will be using the term 'ascending dimension order' when we list the elements of an arrangement of subvarieties of a fixed variety. By this we will understand that if $X_1\subset X_2$ as sets then $X_1<X_2$ i.e. it is a partial order condition with respect to inclusion. We will also often write $N$ in place of the set $\{1,2\dots,n\}$.

\textbf{Acknowledgements}: I am grateful to my advisor, Dan Abramovich, for suggesting this project to me, for numerous conversations and for his constant encouragement to complete it. I would also like to express my gratitude to Angela Gibney and Danny Krashen for their assistance and valuable advice. Finally, I would like to thank Brendan Hassett, Noah Giansiracusa, Qile Chen, Dawei Chen and Sam Molcho for discussing this work and for providing useful comments.

\end{section}

\begin{section}{ A weighted stable compactification}
Throughout this section we will understand the term variety as an integral variety over an algebraically closed field $k$.

\begin{subsection}{Building sets and Wonderful Compactifications}
In this section we recall some definitions and results from \cite {Li}.
\vspace{0.1in}

\begin{mydef}A \textbf{simple arrangement} of subvarieties of a nonsingular variety $Y$ is a finite set $\mathcal{S}=\{S_i\}$ of nonsingular closed  subvarieties $S_i$ properly contained in $Y$ satisfying the following conditions:

\renewcommand{\theenumi}{\roman{enumi}}
\begin{enumerate}
  \item $T(S_i \cap S_j)=TS_i|_{S_i \cap S_j} \cap TS_j|_{S_i \cap S_j}$  and
  \item $S_i \cap S_j$ is either equal to some $S_k \in \mathcal{S}$ or empty.
\end{enumerate}
\end{mydef}
If the smooth varieties $S_i$ and $S_j$ satisfy condition(i) above, we say that $S_i$ and $S_j$ intersect $\textbf{cleanly}$.
\vspace{0.1in}
\begin{mydef} Let $\mathcal{S}$ be an arrangement of subvarieties of $Y$. A subset $\mathcal G \subset \mathcal S$ is called a \textbf{building set} of $\mathcal S $ if for all $S\in\mathcal S$, the minimal elements of $\mathcal G$ containing $S$ intersect transversally and their intersection is equal to $S$ (by convention, the condition is satisfied if $S\in \mathcal G$). These minimal elements are called the $\mathcal G$-factors of $S$.\\
A finite set $\mathcal G$ of nonsingular subvarieties of $Y$ is called a building set if the set of all possible intersections forms an arrangement $\mathcal S$ and if, in addition, $\mathcal G$ is a buiding set of $\mathcal S$. In this case, $\mathcal S$ is called the arrangement induced by $\mathcal G$.
\end{mydef}

For example, let $X$ be a nonsingular variety. The set $\mathcal G =\{\Delta_{12}, \Delta_{13}\}$ of diagonals in $X^3$ is a building set whose induced simple arrangement is the set $\{\Delta_{12}, \Delta_{13}, \Delta_{123}\}$. On the other hand, $\{\Delta_{12}, \Delta_{13}, \Delta_{23}\}$ is not a buiding set; the set of all possible intersections is $\{\Delta_{12}, \Delta_{13}, \Delta_{123}\}$, yet $\Delta_{123} $ is not the transversal intersection of $\Delta_{12}, \Delta_{13}$ and $\Delta_{23}$.\\

\begin{mydef}Let $\mathcal G$ be a building set and set $Y^o:=Y \setminus \bigcup _{G\in \mathcal G} G$. The closure of the image of the locally closed (diagonal) embedding \\
\begin{center}
$Y^o\xhookrightarrow\ Y \times \prod \limits_ {G\in\mathcal G} Bl_G Y$

\end {center}
\vspace{0.1in}
is called the wonderful compactification of $Y$ with respect to $\mathcal G$ and is denoted by $Y_\mathcal G$.
 \end{mydef}
\vspace{0.1 in}
In the sequel, we will frequently make use of the following:
\begin{mydef}The \textbf{dominant transform} of a variety under a blowup is either the strict transform, in case the variety is not contained in the center or the  blowup, or the inverse image of the variety via the blowup map if it is contained in the center.
\end{mydef}
The following theorem (\cite{Li}) will be useful for our construction:
\begin{thm}\label{thm1}
Let $Y$ be a nonsingular variety and $\mathcal G=\{G_1,G_2,\dots ,G_N\}$ be a building set of subvarieties of $Y$. Then,

\begin{enumerate}
\item {the wonderful compactification $Y_\mathcal G$ is a nonsingular variety. Moreover, for each $G\in \mathcal{G}$ there is a nonsingular divisor $D_G\subset Y_{\mathcal{G}}$}, such that
\begin{enumerate}
\item{The union of the divisors is $Y_{\mathcal{G}}\setminus Y^o$}
\item{Any set of these divisors meet transversally. An intesection of divisors $D_{T_1}\cap\dots \cap D_{T_r}$ is nonempty precisely when $\{T_1, \dots ,T_r\}$ is a nest.}
\end{enumerate}

\item if we arrange $\mathcal G$ in such an order that the first $i$ terms  $G_1,G_2,\dots ,G_i$ form a building set for all $1\leq i \leq N$, then \\

\begin{center}
$Y_\mathcal G= Bl_{{{\widetilde{G}}_N}}\cdots Bl_{{{\widetilde{G}}_2}}Bl_{G_1}Y$ 
\end{center}
\vspace{0.1in}
where $Y_\mathcal G$ is the wonderful compactification of $Y$ with respect to $\mathcal G$ and the $\ \widetilde{} \ $  sign on top of each $G_i$ stands for the iterated dominant transform of the latter in the corresponding blowup.
\item let $I_i$ be the ideal sheaf of $G_i$. Then, the wonderful compactification  $Y_\mathcal{G}$ is equal to the blowup of $Y$ along the ideal sheaf $I_1I_2\dots I_N$,\\
\begin{center}
$Y_\mathcal G\iso Bl_{I_N}\dots Bl_{I_2}Bl_{I_1}Y \iso  Bl_{I_1I_2\dots I_N}Y $
\end{center}
\end{enumerate}
\end{thm}
\textbf{Remark}: The divisors $D_G$ are the iterated dominant transforms of the $G_i$ under the blowup  $Y_\mathcal G \rightarrow Y$.\\

\bigskip
The Fulton-Macpherson compactification \cite {FM} is an example of a wonderful compactification as we see in the following example:
\begin{myex}{The Fulton-Macpherson configuration space}
\end{myex}
Consider the configuration space of $n$ distinct labeled points in a nonsingular variety $X$. It is equal to the complement of all diagonals in $X^n$,\\
\begin{center}
 $F(X,n):=X^n\setminus \bigcup \limits_{|I| \geq2} \Delta _I$

\end{center}
\vspace{0.1in}
In \cite{FM}, the compactification of $F(X,n)$, $X[n]$, is constructed inductively as a sequence of blowups. More specifically, $X[1]=X$ and $X[2]$ is the blowup of $X^2$ along the diagonal $\Delta_{12}$. Then $X[3]$ is constructed from $X[2]\times X$ in three steps: First, the exceptional divisor $D$ of $X[2]$ is embedded as a graph of the map from $D$ to $X$ in $X[2]\times X$; the embedded subvariety maps down to $\Delta_{123}$ in $X^3$. Also note that $X[2]$ can be embedded in $X[2]\times X$ in two ways as a graph of the projections from $X[2]$ to $X$; the embedded subvarieties map down to $\Delta_{13}, \Delta_{23}$ in $X^3$. To get $X[3]$ one then needs to blow up along $D$; this gives (the universal family) $X[2]^+$. Then one blows up along the strict transforms of $X[2]$ (embedded as above) which have become disjoint so the order doesn't matter. To get $X[4]$ one then starts over with $X[3]\times X$ by blowing up a subvariety  which maps down to $\Delta_{1234}$, etc. In general, to get $X[n]$ one blows up along subvarieties corresponding to all diagonals in the following order\\

\begin{center}
$\Delta_{12}; \Delta_{123}; \Delta_{13}, \Delta_{23}; \Delta_{1234}; \Delta_{124}, \Delta_{134}, \Delta_{234}; \Delta_{14},  \Delta_{24},   \Delta_{34};\dots $
\end{center}
\vspace{0.1in}
or equivalently $X[n]$ is obtained as a sequence of blow ups along dominant transforms of the above diagonals in the above order. The above sequence satisfies the hypotheses of Theorem 1; for this fact, see the proof of Theorem \ref{thm4} below. We therefore conclude that $X[n]$ is the wonderful compactification $X^n_{\mathcal G}$ of $X^n$ with respect to the building set $\mathcal G$ of all diagonals. This recovers \cite[Proposition 4.1]{FM}.  \\

Moreover, theorem 1 implies that $X[n]$ can also be obtained by the following symmetric sequence of blowups along dominant transforms of diagonals in ascending dimension\\
\begin{center}
$\Delta_{12\dots n}; \Delta_{12\dots (n-1)},\dots ,\Delta_{23\dots n}; \dots,   \Delta_{12},\dots   ,\Delta_{(n-1)n} $
\end{center}
\vspace{0.1in}
For more details see the proof of Lemma 1 below.
\end{subsection}
\begin{subsection}{Weighted compactification of configuration spaces}\label{sec2.2}In \cite{FM}, the authors describe a natural and beautiful way of compactifying the configuration space of $n$ distinct labeled points in a nonsingular variety $X$. The degenerate objects of the resulting compactification consist of varieties with  $n$ points that remain distinct. In this section we introduce an alternate compactification of this configuration space, by equipping points with weights and allowing them to collide in the degenerations depending on the accumulation of their weights. 

\begin{subsubsection}{n-pointed $\mathcal{A}$ stable degenerations} We now give a geometric description of the degenerate objects of our compactification(see also \cite{FM},\cite{Pa1}). 
Let $X$ be a nonsingular variety over an algebraically closed field $k$ and $\mathcal{A}:=\{a_1,a_2, \dots,a_n\} $ be a set of rational numbers such that\\
\renewcommand{\labelitemi}{$\bullet$}
\begin{itemize}
\item $0<a_i\leq 1 , i=1,2,\dots,n$ 
\end {itemize}
\vspace{0.1in}
\vspace{0.1in}
 Let $(x_1,x_2,\dots, x_n)$ be an ordered $n$-tuple of points of $X$ carrying weights $a_i$. A subset $S\subset N$ is said to be \textit{coincident} at $x\in X$ if 
 \begin{itemize}
 \item $|S|\geq2$,
 \item $\sum \limits _{i\in S} a_i>1$ and
 \item for all $i\in S$, $x_i=x$
 \end{itemize} 
 \vspace{0.1in}
 A \textit{screen} of $S$ at $x$ consists of the data $(t_i)_{i\in S}$ such that:
 \begin{enumerate}
\item $t_i\in T_x$, the tangent space of $X$ at $x$
\item there exist $i,j\in S$ such that $t_i\neq t_j$ 
\end{enumerate}
Two data sets $(t_i)_{i\in S}$ and $(t_i')_{i\in S}$ are equivalent if there exist $c\in\CC^*$ and $v\in T_x$ s.t. $c\cdot t_i +v =t_i'$ for all $i\in S$. Now, consider the $n$-tuple $(x_1,x_2,\dots x_n)$ together with the collection of all coincident sets $S$. We construct the (weighted) $n$-pointed $\mathcal{A}$ stable degeneration of $X$ as follows. Let $z$ be a coordinate that occurs multiple times in $(x_1,\dots x_n)$ and generates a coincident set at $z$. We blow up $X$ at $z$ and attach the projective completion $\PP(T_z\oplus1)$ along the exceptional divisor $\PP(T_z)$, which is identified with the infinity section. Note that the complement $\PP(T_z\oplus1)\setminus \PP(T_z)$ is isomorphic to the affine space $T_z$. Let $S_z$ be the maximal coincident subset at $z$. The screen corresponding to $S_z$ associates points of $T_z$ to the indices that lie in $S_z$. By condition (2) for screens, we see that some separation of those points occurs. We continue this process by blowing up points in the new spaces $T_z$ specified by the subsequent screens until all screens have been used, for all such coordinates $z$. The resulting variety is equipped with $n$ points $s_i$ lying in the smooth locus. By this description we see that if $S\subset N$ and $\sum\limits_{i\in S}a_i>1$ then some separation of the sections $(s_i)_{i\in S}$ necessarily occurs. This means that if the sections $(s_i)_{i\in S}$ all coincide for some $S$, then $\sum\limits_{i\in S}a_i\leq1$.\\

To any $\mathcal{A}$ stable degeneration we associate a tree whose vertices correspond one to one to its components and whose vertices correspond one to one to the intersections of its components. There exists a distinguished component called the \textit{root}, which is a blowup of $X$ at a finite set of points. We call a component an \textit{end} if it is not the root and the valence of its vertex is equal to 1. A \text{ruled} component is one whose vertex has valence 2. End components come with at least three distinct markings: at least two coming from distinct smooth sections and one (node) from an intersection with another component. Ruled components also come with three distinct markings: at least one from a smooth section and two from intersections with other components (nodes). In other words all components, except for the root, come with at least three distinct markings, so there are no nontrivial automorphisms of an $n$- pointed $\mathcal{A}$ stable degeneration preserving the root. This justifies the term \textit{stability}.
\end{subsubsection}
\begin{subsubsection}{The weighted compactification}
\bigskip
Let \\

\begin{center}
$\mathcal{G}:=\{ \Delta_{12\dots n}; \Delta_{12\dots (n-1)},\dots ,\Delta_{23\dots n}; \dots;  \Delta_{12},\dots   ,\Delta_{(n-1)n} \} $
\end{center}
\bigskip
be the set of diagonals of $X^n$ listed in ascending dimension order.
\vspace{0.1in}

Also, let 
\begin{center}
$\mathcal{G}_{\mathcal{A}}:= \{ \Delta_S\subset X^n | S\subset N \,\text{and}\, \sum \limits _{i\in S} a_i >1\}$

\end{center}
 and list its elements in ascending dimension order. \\
\begin{lem}\label{lem1} $\mathcal{G}_{\mathcal{A}}$ satisfies the condition of Theorem 1 (2). In particular, $\wbuild$ is a building set.

\end{lem}
\textbf{Proof}: Consider the set of the first $i$ elements,
${\mathcal{G}_{\mathcal{A}}}_{,i}$. Pick an arbitrary subset of $k$ elements of ${\mathcal{G}_{\mathcal{A}}}_{,i}$, say $\{\Delta_{I_1}, \Delta_{I_2},\dots , \Delta_{I_k}\}$  and consider the intersection $S_k:=\Delta_{I_1}\cap \Delta_{I_2}\cap\dots \cap \Delta_{I_k}$, i.e. an element of the simple arrangement induced by ${\mathcal{G}_{\mathcal{A}}}_{,i}$. We need to show that the minimal elements of ${\mathcal{G}_{\mathcal{A}}}_{,i}$ that contain $S_k$ intersect transversally along $S_k$. Observe that the above intersection can be written uniquely as $\Delta_{J_1}\cap\Delta_{J_2}\dots\cap \Delta_{J_{k'}}$, where $k'\leq k$ and for  $j=1,\dots k'$, the indices $J_j$ are unions of some of the $I_l$ and also pairwise disjoint. Indeed, we start with $I_1$; if $I_1$ is disjoint from the union $I_2\cup \dots \cup I_k$, we set $J_1:=I_1$. Otherwise, we pick the first element (index) in the collection of indices $\{I_1, I_2,\dots I_k\}\setminus\{I_1\}=\{I_2,\dots I_k\}$ that is not disjoint (as a set) from $I_1$ and we call this $I_{t_1}$. Then we pick the first element in the collection $\{I_1, I_{2},\dots I_k\}\setminus \{I_1, I_{t_1}\}$ that is not disjoint (as a set) from the union $I_1\cup I_{t_1}$ and we carry on until we have chosen $I_1, I_{t_1}, \dots I_{t_l}$ such that all elements of the set $\{I_1, I_{2},\dots I_k\}\setminus \{I_1, I_{t_1}, \dots I_{t_l}\}$ are disjoint from the union $I_1\cup I_{t_1} \cup \dots \cup I_{t_l}$. We set $J_1:= I_1\cup I_{t_1} \cup \dots \cup I_{t_l}$.  Of course, $\Delta _{J_1}= \Delta_ {I_1} \cap \Delta_{I_{t_1}} \dots \cap \Delta_{I_{t_l}}$. Next, to obtain $J_2$, we repeat the same procedure for the first index among $\{I_1, I_{2},\dots I_k\}\setminus \{I_1, I_{t_1}, \dots I_{t_l}\}$, etc. and we continue until all elements of the set $\{I_1, I_2, \dots I_k\}$ have been exhausted. At the end of the procedure we will obtain indices $J_i, \, i=1,\dots k'$, where $k'\leq k$ with the property that $S_k:=\Delta_{I_1}\cap\Delta_{I_2}\dots\cap \Delta_{I_{k}}=\Delta_{J_1}\cap\Delta_{J_2}\dots\cap \Delta_{J_{k'}}$. Since the $J_i$ are pairwise disjoint, it is clear that this expression is unique.  \\

Observe that $\Delta_{J_j} \in {\wbuild}_{,i} $ for all $j=1,\dots k'$: indeed, for each $j$, the index set ${J_j}$ contains some $I_l, l\in\{1,2,\dots k\}$, so \\
\begin{center}
 $\sum \limits _{i\in J_j} a_i \geq \sum \limits _{i\in I_l} a_i>1$
\end{center}
\vspace{0.1in}
Therefore, the minimal elements of ${\wbuild}_{,i}$ that contain $S_k=\Delta_{J_1}\cap\Delta_{J_2}\dots\cap \Delta_{J_{k'}}$ are precisely the elements $\Delta_{J_1}, \Delta_{J_2},\dots \Delta_{J_{k'}}$, which are seen to intersect transversally, since the $J_i$ are disjoint.\qed

\vspace{0.1in}
\begin{mydef} The weighted configuration space with respect to $\mathcal{A}$ is the complement of the union of the diagonals that belong to $\mathcal{G}_{\mathcal{A}}$ in $X^n$\\
 \begin{center}
 $F_\mathcal{A}(X,n):= X^n\setminus \bigcup\limits_ {\Delta_I \in \mathcal{G}_{\mathcal{A}}} \Delta_I$ . \\
  \end{center}
\end{mydef}
 
  \vspace{0.1 in}

The weighted configuration space may be viewed as the parameter space of $n$ labeled points in $X$ carrying weights $a_i$ subject to the following condition:\\
\renewcommand{\labelitemi}{$\bullet$}
\begin{itemize}
\item for any set of labels $S\subset N$ of coincident points we have $\sum \limits _{i\in S} a_i \leq 1\ $.
\end {itemize}
 
 We are now ready to define its weighted compactification:
 
 \begin{mydef}The weighted compactification of the configuration space $F_\mathcal{A}(X,n)$ is the wonderful compactification of $X^n$ with respect to the building set $\mathcal{G}_{\mathcal{A}}$, that is the closure of the image of the diagonal embedding\\
 
  \begin{center}
 
$ X^n\setminus \bigcup\limits_ {\Delta_I \in \mathcal{G}_{\mathcal{A}}} \Delta_I\xhookrightarrow\ X^n \times\displaystyle \prod \limits_ {\Delta_I\in\ \mathcal{G}_{\mathcal{A}}} Bl_{\Delta_I} X^n$

\end {center}
\vspace{0.1in}
We denote the weighted compactification by $X_\mathcal{A}[n]$.
\end{mydef}
\vspace{0.1in}
\textbf{Remark}: $\wbuild$ is empty if $\sum \limits_{i=1}^{n} a_i\leq1$. In this case $\wmod$ is trivially equal to $X^n$.\\

\textbf{Remark}: The variety $\wmod$ is also a compactification of the configuration space of $n$ \textit{distinct} labeled points $F(X,n)=X^n\setminus \bigcup\limits_ {|I|\geq2} \Delta_I$, since the latter is an open subset of 
$X^n\setminus \bigcup\limits_ {\Delta_I \in \mathcal{G}_{\mathcal{A}}} \Delta_I$, so the closure of \\
\begin{gather*}
X^n\setminus \bigcup\limits_ {|I|\geq2} \Delta_I\xhookrightarrow\ X^n\setminus \bigcup\limits_ {\Delta_I \in \mathcal{G}_{\mathcal{A}}} \Delta_I \xhookrightarrow\ 
X^n \times\displaystyle \prod \limits_ {\Delta_I\in\ \mathcal{G}_{\mathcal{A}}} Bl_{\Delta_I} X^n
\end{gather*}
is the same as the closure of \\
 
  \begin{center}
 
$ X^n\setminus \bigcup\limits_ {\Delta_I \in \mathcal{G}_{\mathcal{A}}} \Delta_I\xhookrightarrow\ X^n \times\displaystyle \prod \limits_ {\Delta_I\in\ \mathcal{G}_{\mathcal{A}}} Bl_{\Delta_I} X^n$

\end {center}

For this reason we can consider it a generalization of the original Fulton-Macpherson compactification. \\

The variety $\wmod$ enjoys the following nice properties:
 \begin{thm}
\begin{enumerate}
\item {$\wmod$ is a nonsingular variety. The boundary  $X_\mathcal{A}[n] \setminus (X^n \setminus \bigcup \limits_ {\Delta_I \in \mathcal{G}_{\mathcal{A}}} \Delta_I) $ is the union of $|\mathcal{G}_{\mathcal{A}}|$ divisors $D_I,$ where $I\subset N   , |I|\geq2$ and $ \sum \limits _{i\in I} a_i >1\ $.}

\item{ Any set of boundary divisors intersects transversally. An intersection of divisors $D_{I_1}\cap D_{I_2}\cap\dots D_{I_k}$ is nonempty precisely when the sets are nested in the sense that any pair $\{I_i,I_j\}$ either has empty intersection or one set is contained in the other. }

 \item{$X_\mathcal{A}[n]$ is the iterated blowup of $X^n$ along dominant transforms of elements of $\mathcal{G}_{\mathcal{A}}$}.
\end{enumerate}
\end{thm} 

 \textbf{Proof}: This is an immediate consequence of Theorem \ref{thm1} and Lemma \ref{lem1}. In part (3), in particular,  the order of the blowups agrees with the order we wrote down the elements of $\wbuild$.  
 \qed
\vspace{0.1in}

The following basic lemma will be useful for the proofs of the results that appear in the sequel.
\begin{lem}\label{lem2}Let  $Z$ be a smooth subvariety of a smooth variety $Y$ and let $\pi: Bl_Z Y\rightarrow Y$ be the blowup, with exceptional divisor $E=\pi^{-1}(Z)$.\\ 
\begin{enumerate}
\item Let $V$ be a smooth subvariety of $Y$, not contained in $Z$, and let $\widetilde{V}\subset Bl_Z Y$ be its strict transform. Then,\\
\begin{enumerate}
\item if $V$ meets $Z$ transversally (or is disjoint from $Z$), then $\widetilde{V}=\pi^{-1}(V)$ maps isomorphically to $V$ and $I_{\pi^{-1}(V)}=I_{\widetilde{V}}$.\\
\item if $V\supset Z$, then $I_{\pi^{-1}(V)}=I_{\widetilde{V}}\cdot I_E$. If Z has codimension 1 in V, then the projection from $\widetilde{V}$ to $V$ is an isomorphism. \\
\end{enumerate}
\item If $Z_i$ are smooth subvarieties of $Y$ such that $\sum I_{Z_i}= I_Z $, where $I_{Z_i}$ (resp. $I_Z$) are the ideal sheaves of $Z_i$ (resp. $Z$), then the strict transforms of the $Z_i$ in $Bl_Z Y$ do not intersect simultaneously.  \\
\item Let $Z_1, Z_2$ be smooth subvarieties of $Y$ intersecting transversally. \\
\begin{enumerate}
\item Assume $Z_1\cap Z_2 \supsetneq Z$. Then their strict transforms $\widetilde{Z_1}$ and $\widetilde{Z_2}$ intersect transversally and $\widetilde{Z_1}\cap \widetilde{Z_2}= \widetilde{Z_1\cap Z_2}$.\\
\item Assume $Z$ intersects transversally with $Z_1$ and $Z_2$, as well as with their intersection $Z_1\cap Z_2$. Then their strict transforms $\widetilde{Z_1}$ and $\widetilde{Z_2}$ intersect transversally and $\widetilde{Z_1}\cap \widetilde{Z_2}= \widetilde{Z_1\cap Z_2}$.\\
\item If $Z_1\supseteq Z$ and $Z_2$ intersects transversally with $Z$, then the intersections \\
\begin{gather*}
\widetilde{Z_1}\cap \widetilde{Z_2}\,\, \text{and}\,\, (E\cap \widetilde{Z_1})\cap \widetilde{Z_2}
\end{gather*}  
are transversal. Moreover, $\widetilde{Z_1}\cap \widetilde{Z_2}=\widetilde{Z_1\cap Z_2}$.\\
\end{enumerate}
\item Let $W$ be a smooth subvariety of $Y$ that intersects transversally with $Z$. Then,
\begin{gather*}Bl_{\widetilde{W}} Bl_Z Y=Bl_{\widetilde{Z}} Bl_W Y\\
\end{gather*}
\end{enumerate}
\end{lem}
\textbf{Proof}: (1) is the content of \cite[Lemma 4.1(a),(b)]{FM}, (2) and (4) are standard, while (3) is essentially the content of \cite[Lemma 2.9,(iii), (iv) and (v)]{Li}. 
\qed
\vspace{0.1in}

\begin{prop}\label{prop1}Let $\pi:\wmod \rightarrow X^S$ be the composition of the blowup map (see theorem 2(3)) $\wmod \rightarrow X^n$ with the projection $X^n\rightarrow X^S$ onto the coordinates labeled by $S$ for some $S\subset N$ such that $\Delta_S\in \wbuild$. Then, 
\begin{center}
$I_{{\pi^{-1}}\Delta_S}=\displaystyle\prod\limits_{T\supset S} I_{D_T}$
\end{center}
where $\Delta_S$ is the small diagonal in $X^S$.
\end{prop}
\textbf{Proof}:We extend the partial order of the set of indices of the elements of $\wbuild$ to a total order which respects inclusions of diagonals:
\begin{center}
 $I_0:=(12\dots n)<I_1<\dots <I_{|\wbuild|-1}$ \\
\end{center}
Define $X_0:=Bl_{{\Delta_{I_0}}}X^n$ and let $\wbuild^{(I_0)}:=\{ {\Delta_{I_0}}^{(I_0)} , {\Delta_{I_1}}^{(I_0)}, \dots {\Delta_{I_l}}^{(I_0)}, \dots ,{\Delta_{I_k}}^{(I_0)},\dots\}$ be the set of dominant transforms of the varieties(diagonals) that belong to $\wbuild$ in $X_0$. Then define $X_1:=Bl_{{\Delta_{I_1}}^{(I_{0})}}X_0$ and let  $\wbuild^{(I_1)}:=\{ {\Delta_{I_0}}^{(I_1)} , {\Delta_{I_1}}^{(I_1)}  \dots {\Delta_{I_l}}^{(I_1)}, \dots ,{\Delta_{I_k}}^{(I_1)},\dots\}$ be the set of dominant transforms of the varieties that belong to $\wbuild^{(I_0)}$ in $X_1$. For $k\geq2$, we may define inductively  the $I_k$-th iterated blowup, \\
\begin{center}
$X_k:= Bl_{\Delta_{I_{k}}^{(I_{k-1})}} X_{k-1}= Bl_{{\Delta_{I_k}}^{(I_{k-1})}}\dots Bl_{{\Delta_{I_0}}}X^n$
\end{center}
\vspace{0.1in}
and $\wbuild^{(I_k)}:=\{ {\Delta_{I_0}}^{(I_k)} ;\dots {\Delta_{I_l}}^{(I_k)}, \dots ,{\Delta_{I_k}}^{(I_k)},\dots\}$ as the set of dominant transforms of the subvarieties of $X_{k-1}$ that are the elements of $\wbuild^{(I_{k-1})}$. In other words, $\wbuild^{(I_k)}:=\{ {\Delta_{I_0}}^{(I_k)} ,\dots {\Delta_{I_l}}^{(I_k)}, \dots ,{\Delta_{I_k}}^{(I_k)},\dots\}$ is the set of the $I_k$-th iterated dominant transforms of the diagonals of $\wbuild$ in $X_k$.\\
\vspace{0.1in}
Then, by \cite[Remark after Definition 2.12]{Li}, we deduce that  each $\wbuild^{(I_k)}$ is a building set.\\
Let $k$ be such that $I_k=S$. We prove the following claim by induction :\\

$\mathit{Claim}$: Let $\pi_{j+1,j}:X_{j+1}\rightarrow X_j$ and  $\pi_j: X_j\rightarrow X^n$ be the natural blowup maps. Then, for any $j$, 
\begin{center}
$\displaystyle I_{{\pi_j}^{-1}({\Delta_{I_k}})}= \displaystyle I_{{\Delta_{I_k}}^{(I_{j})}} \cdot  \prod\limits_ {\substack{ l\leq j\\ I_l\supsetneq I_k}} I_{{\Delta_{I_l}}^{(I_j)}}$
\end{center}
\vspace{0.1in}
$\mathit{Proof}$: Let $j=0$. If $k=0$ there's nothing to prove. Otherwise, $\Delta_{I_k}\supsetneq \Delta_{I_0}$, so by Lemma 2, we have
\begin{center}
$\displaystyle I_{{\pi_0}^{-1}({\Delta_{I_k}})}= I_{{\Delta_{I_k}}^{(I_{0})}}\cdot  I_{{\Delta_{I_0}}^{(I_{0})}} $
\end{center}

\bigskip
Now assume the claim is true for $j$. By pulling back to $X_{j+1}$ via $\pi_{j+1,j}:X_{j+1}\rightarrow X_{j}$, we obtain \\
\begin{gather}
I_{{\pi_{j+1}}^{-1}({\Delta_{I_k}})}=I_{{{\pi_{j+1,j}}^{-1}}({\Delta_{I_k}}^{(I_j)})}\cdot \displaystyle\prod\limits_ {\substack{ l\leq j\\ I_l\supsetneq I_k}} I_{{{\pi_{j+1,j}}^{-1}}({\Delta_{I_l}}^{(I_j)})}\end{gather}
Observe that each of the terms ${\Delta_{I_l}}^{(I_j)}$ are divisors so they cannot be contained in the center ${\Delta_{I_{j+1}}}^{(I_j)}$ due to dimension reasons (note $l<j+1$). Therefore, the center is minimal in $\wbuild^{(I_{j})}$, so by  \cite[Lemma 2.6(i)]{Li}, the above terms must be  transversal to the center, hence, by Lemma \ref{lem2}(1) \\
\begin{align}
I_{{{\pi_{j+1,j}}^{-1}}({\Delta_{I_l}}^{(I_j)})}=I_{{\Delta_{I_l}}^{(I_{j+1})}}
\end{align}

\bigskip
Moreover, again by using Lemma 2 we get
\begin{align}
I_{{{\pi_{j+1,j}}^{-1}}({\Delta_{I_k}}^{(I_j)})}=
\begin{cases}
I_{{\Delta_{I_k}}^{(I_{j+1})}}\cdot I_{{\Delta_{I_{j+1}}}^{(I_{j+1})}}  &\text{, if}\, I_{j+1} \supsetneq I_k\\
I_{{\Delta_{I_k}}^{(I_{j+1})}} &\text{, otherwise}
\end{cases}
\end{align}

\bigskip

Plugging (2) and (3) into (1), we obtain\\
\begin{align*}
I_{{\pi_{j+1}}^{-1}({\Delta_{I_k}})}=
\begin{cases}
I_{{\Delta_{I_k}}^{(I_{j+1})}}\cdot I_{{\Delta_{I_{j+1}}}^{(I_{j+1})}} \cdot \displaystyle\prod\limits_ {\substack{ l\leq j\\ I_l\supsetneq I_k}} I_{{\Delta_{I_l}}^{(I_{j+1})}} &\text{, if}\, I_{j+1} \supsetneq I_k\\[3em]
I_{{\Delta_{I_k}}^{(I_{j+1})}}  \cdot \displaystyle\prod\limits_ {\substack{ l\leq j\\ I_l\supsetneq I_k}} I_{{\Delta_{I_l}}^{(I_{j+1})}}&\text{, otherwise}
\end{cases}
\end{align*}

\bigskip

In both cases, we have $I_{{\pi_{j+1}}^{-1}({\Delta_{I_k}})}=I_{{\Delta_{I_k}}^{(I_{j+1})}}  \cdot \displaystyle\prod\limits_ {\substack{ l\leq j+1\\ I_l\supsetneq I_k}} I_{{\Delta_{I_l}}^{(I_{j+1})}}$, so
 we conclude the claim for $j+1$.\\

By the remark after Theorem 1, the divisors $D_T$ in $\wmod$ are the iterated dominant transforms of $\Delta_T\in \wbuild$ under the sequence of blowups $\wmod\rightarrow \ X^n$. By taking $j$ such that $I_j$ is the greatest index (with the above total order) in the set of indices of the elements of $\wbuild$, in other words  $I_j=I_{|\wbuild|-1}$, and applying the claim we conclude the proposition.\qed

\vspace{0.1in}
\textbf{Remark}: It does not matter which total order we choose, as long as it respects the dimension order of the diagonals; in this case $\wbuild$ satisfies the condition of Theorem \ref{thm1}(2) (as in lemma 1), so the resulting blowup is independent of the ordering.\\

\textbf{Remark}: The above Proposition is the analogue of \cite[Theorem 3(3)]{FM}. \\

\textbf{Remark}: In case $a_i=1$ for all $i=1,\dots,n$, the moduli space  $X_\mathcal{A}[n]$ coincides with the Fulton Macpherson  compactification. In general, the weighted moduli space enjoys more properties analogous to those of  the original Fulton Macpherson space \cite{FM}, as we shall see in the following sections.
\end{subsubsection}
\end{subsection}
\begin{subsection}{The universal family}\label{sec2.3}The content of this section is summarized in the following Theorem:

\begin{thm}\label{thm3}
 There exists a `universal' family  $X_\mathcal{A}[n]^+\rightarrow X_\mathcal{A}[n]$ equipped with $n$ sections 
 $\sigma_i : X_\mathcal{A}[n]\rightarrow X_\mathcal{A}[n]^+, i=1,\dots n$. It is a flat morphism between nonsingular varieties, whose fibers are n pointed $\mathcal{A}$ stable degenerations of $X$.
\end{thm}
\vspace{0.1in}

The remainder of this section is devoted to the proof of the above theorem. We first construct the universal family as a sequence of blowups and define auxiliary spaces. The proof is broken into various steps. We also prove various facts about the geometry of the intersections of the subvarieties involved in each blowup.
\vspace{0.1in}

Consider the product $Y_0:=\wmod \times X$. For any $I\subset N$, let $I_+\subset N_+:=\{1,2,\dots,n+1\}$ be the set $I\cup \{n+1\}$. We define various subvarieties of $Y_0$. First, consider the divisors $D_I$ of $\wmod$, i.e. the iterated dominant transforms of the elements of $\wbuild$ under the sequence of  blowups $\wmod\rightarrow X^n$. Those map to $\Delta_I$ via $\wmod \rightarrow X^n\rightarrow X^I$, where the second map is the projection onto the coordinates labeled by $I$. 
\begin{itemize}
\item With $\Delta_I$ identified with $X\subset X^I$, we embed each $D_I$ in $Y_0:=\wmod\times X$ as the graph  of the map $D_I\rightarrow X$. The image of this embedding is defined to be $D_{I+}^{(Y_0)}$.
\item We define ${D_I}^{(Y_0)}:=D_I\times X$.
\item For each $i\in N$ we define $\Delta_{i+}^{(Y_0)}$ to be the embedding of $\wmod$ in $\wmod\times X$ given by the graph of the morphism  $\wmod\rightarrow X^n\xrightarrow{p_i} X$, where $p_i$ is the projection to the $i$-th factor.\end{itemize}
 
 \vspace{0.1in}
Define $Y_1$ to be the blowup of $Y_0$ along $D_{N+}^{(Y_0)}$. Assume $Y_k$ is defined. We will call the variety $Y_k$ the \textbf{$k$-th step}. Then we define $Y_{k+1}$ to be the blowup of $Y_k$ along the disjoint union of the subvarieties $D_{I+}^{(Y_k)}$, where $|I|=n-k$ and $\Delta_I \in \wbuild$; those subvarieties will be shown to be disjoint in the following Propositions.  For any $I, I_+\subset N_+$ s.t $\Delta_I\in\wbuild$
we define the subvarieties $D_{I}^{(Y_{k+1})}$ and $D_{I+}^{(Y_{k+1})}$ of $Y_{k+1}$  to be the dominant transforms of the subvarieties $D_{I}^{(Y_{k})}$ and $D_{I+}^{(Y_{k})}$ respectively of $Y_k$. We do the same for $\Delta_{i+}^{(Y_{k+1})}$.  Let $s$ be the greatest integer for which there exist a set $I\subset N$ such that $|I|=n-s$ and $\Delta_I\in \wbuild$. Then, the \textbf{universal family} is defined as\\
\begin{center}
$\wmod^+:=Y_{s+1}$
\end{center}
\vspace{0.1in} 
In other words, if we define the set of subvarieties  ${\wbuild}_{+}:=\{D_{I+}\subset \wmod\times X \, |\Delta_I \in \wbuild\}$ and list its elements in order preserving bijection with the elements of  $\wbuild$, that is, $D_{I+}<D_{J+}$ if and only if $\Delta_I<\Delta_J$, then $\wmod^+$ is the iterated blowup along dominant transforms of the elements of ${\wbuild}_{+}$ in that order (note that ${\wbuild}_{+}$ is \textbf{not} a building set). \\

Finally, we define  \textbf{stage} $(k,S)$ to be the blowup of $\wmod\times X$ along the iterated dominant transforms of the subset of ${\wbuild}_{+}$ obtained by cutting out all elements listed after $D_{S+}$, where $|S|=n-k$. The $Y_k$ defined above are special cases of the $(k,S)$ stages; for instance, $Y_1$ is stage $(0,N)$. The dominant transforms of the initial subvarieties $D_I^{(Y_0)}, D_{I+}^{(Y_0)}$ and ${\Delta}_{i+}^{(Y_0)}$ of $\wmod\times X$ in the $(k,S)$ stage are denoted by $D_I^{(S+)}, D_{I+}^{(S+)}$ and $\Delta_{i+}^{(S+)}$ respectively.  \\

A few propositions about the geometry of intersections of the subvarieties of $\wmod^+$ as well as the intermediate steps $(k,S)$ are in order. First, let us recall the notation of the proof of Proposition \ref{prop1}:\\

We extend the partial order of the set of indices $\{(12\dots n);(12\dots n-1),\dots (23\dots n); \dots \}$ of the elements of $\wbuild$ to a total order:
\begin{center}
 $I_0:=(12\dots n)<I_1:=(12\dots n-1)<I_2:=(23\dots n)<\dots <I_{|\wbuild|}$ \\
\end{center}
which respects inclusions of diagonals.
\vspace{0.1in}
Based on this, we give the following :\\
\begin{mydef}A \textbf{nest} in $N_+$ is a collection $\mathcal{S}$ of subsets of $N_+$ such that each pair of subsets is either disjoint or one contains the other and for any $S\in\mathcal {S}$, we have $\Delta_{S\cap N}\in \wbuild$. Let $S\subset N$ such that $\Delta_S\in \wbuild$. We say that a nest $\mathcal{S}$ satisfies condition $(*)_S$ if:\\
\begin{itemize}
\item[$(*)_S$]
 $\mathcal{S}$ contains at most one set of the form $I_+$ with $I>S$ with respect to the above order and if it contains such a set $I_+$, then it does not contain $I$. \\ 
 \end{itemize}
 \end{mydef}
 Now, let $S\subset N$ such that $|S|=n-k$ and assume that ${\{D_J^{(S+)}\}_{J\in\mathcal{S}}}$ is a collection of subvarieties of stage $(k,S)$. We have the following:

\begin{prop}\label{prop2}
Assume that $\mathcal{S}$ satisfies $(*)_S$. Then, 
 the subvarieties $D_J^{(S+)}$ of stage $(k,S)$, where $J\in\mathcal{S}$, intersect transversally.
  \end{prop}
\bigskip
\begin{prop}\label{prop3}Consider $I,S\subset N$ such that $\Delta_I,\Delta_S \in \wbuild$ and $I<S$. Then, for any $i\in I$, the subvarieties  $D_I^{(S+)}$ and $\Delta_{i+}^{(S+)}$ of stage $(k,S)$ intersect transversally along $D_{I+}^{(S+)}$: \\
\begin{gather*} D_I^{(S+)}\cap \Delta_{i+}^{(S+)}=D_{I+}^{(S+)}
\end{gather*}    
\end{prop}
The above propositions are similar to the \cite[Proposition 3.1, Proposition 3.3(b)]{FM} and their proofs follow from the definitions and Lemma \ref{lem2}. We also refer the reader to Lemma \ref{lem8} below for a similar proof. \\
The crucial step towards for the proof of Theorem \ref{thm3} is the following:

\bigskip

\begin{prop}\label{prop4}
\renewcommand{\theenumi}{\roman{enumi}}
\begin{enumerate}
\item (stability): At stage $(k,S)$ we have: 
$\displaystyle \bigcap\limits _{i\in S, |S|=n-k}\Delta_{i+}^{(S+)}=\emptyset$ .
\item(separation of centers for $Y_{k+1}$): $D_{S+}^{(Y_{k+1})}\cap D_{T+}^{(Y_{k+1})}=\emptyset$ whenever $|S|=|T|=n-k-1, S\neq T$ and $\Delta_S,\Delta_T \in \wbuild$.\\

\item $D_{S+}^{(Y_{k+1})}\cap D_{T+}^{(Y_{k+1})}=\emptyset$ whenever neither $S$ nor $T$ are contained in each other.
\end{enumerate}
\end{prop}
\vspace{0.1in}
\textbf{Proof}:  By the definitions and proposition \ref{prop1}, we have:\\
\begin{gather*}
\displaystyle \bigcap\limits _{i\in S, |S|=n-k}\Delta_{i+}^{(Y_0)}=\bigcup\limits_{N\supset T \supset S}D_{T+}^{(Y_0)}
 \end{gather*}
 
 Note that we may prove (i) over open covers of $\wmod\times X$. To choose an appropriate open cover, we define a \textbf{chain} between $N$ and $S$ to be a collection of sets:\\
  \begin{gather*}
N:=S_0\supset S_1 \supset \dots \supset S_k:=S
\end{gather*}
such that $|S_i|=n-i$ for $i=0,\dots,k$. Let $\ch$ be the set of all chains between $N$ and $S$. For every chain $C\in\ch$, define the \textbf{complement} of $C$ to be the set $C':=\{T\,| N\supset T\supset S \,\text{and}\,T\notin C\}$ and denote by $\cc$ the collection of all such $C'$. Assume $|\cc|=m$ and let $\cc:=\{C'_1, C'_2,\dots C'_m\}$ be the set of all complements of the chains in $\ch$.  Now, let
\begin{gather*}
 U_C:=\wmod\times X \setminus \bigcup\limits _{I\in C'} D_{I+}^{(Y_0)}
\end{gather*}
The collection of sets $\{U_C\}_{C\in\ch}$ is an open cover of $\wmod \times X$. Indeed, \\
\begin{gather*}
\bigcup\limits_{C\in \ch}U_C=
\bigcup\limits_{C'\in \cc}\left(\wmod\times X \setminus \bigcup \limits_{I\in C'} D_{I+}^{(Y_0)}\right)\\
=\wmod\times X\displaystyle \setminus \bigcap \limits_{C'\in \cc} \left(\bigcup\limits _{I\in C'} D_{I+}^{(Y_0)}\right)\\
=\wmod\times X\displaystyle \setminus \bigcup \limits_{\substack{\{I_i\}_{i=1}^{m} \in \prod\limits _{i=1}^{m} C'_i  }} \left(\bigcap\limits_{i=1}^{m}{D_{{I_i}+}}^{(Y_0)}\right)
\end{gather*}
 Assume that there exists an intersection $\bigcap\limits_{i=1}^{m}{D_{{I_i}+}}^{(Y_0)}$ that is nonempty. By Theorem 2(2), we deduce that such an intersection is nonempty precisely when the $I_i$ are nested in the sense that each pair $\{I_i,I_j\}$ either has empty intersection or one set is contained in the other. Since all $I_i$ above contain $S$, they can't be disjoint, so their intersection can be nonempty only if the indices $I_i$ appearing above are part of a chain between $N$ and $S$, which we denote by $C_{i_0}$. However, in this case the complement $C_{i_0}'$ of $C_{i_0}$ appears as a factor in the product $\prod\limits _{i=1}^{m} C'_i $, so there must be an index $I_{i_0}\in\{I_1,\dots I_m\}$ such that $I_{i_0}\in C_{i_0}'$, which is a contradiction.  \\

Restricting to the above open cover allows us to prove the disjointness of (i) after blowing up at a sequence of iterated strict transforms of ${D_{I+}}^{(Y_0)}$ indexed by a chain between $N$ and $S$. More specifically, let \\
  \begin{gather*}
N:=S_0\supset S_1 \supset \dots \supset S_k:=S
\end{gather*}
be a chain and set $D_{S_i}^{(0)}:={D_{S_i}}^{(Y_0)}$, $D_{S_{i+}}^{(0)}:={D_{S_{i+}}}^{(Y_0)}$ and $\Delta_{i+}^{(0)}=\Delta_{i+}^{(Y_0)}$. Now consider the ordered subset of $\wmod\times X$: \\
\begin{gather*}
\{D_{S_{0+}}^{(0)}, D_{S_{1+}}^{(0)}, \dots , D_{S_{k+}}^{(0)}\}
\end{gather*}
We successively perform blow ups at the iterated dominant transforms of the above set in the above order: let $X_0=\wmod\times X$. Set $X_1=Bl_{D_{S_{0+}}^{(0)}} X_0$ and define $D_{S_0}^{(1)}, D_{S_{0+}}^{(1)}$ and  $\Delta_{i+}^{(1)}$ to be the dominant transforms of $D_{S_0}^{(0)}, D_{S_{0+}}^{(0)}$ and  $\Delta_{i+}^{(0)}$ respectively (for $i\in S$). Inductively we define $X_{j+1}=Bl_{D_{S_{j+}}^{(j)}} X_j$  and denote by $\pi_{j+1}:X_{j+1}\rightarrow X_j$ the blowup map. For each  $i,j=0,\dots k$, let $D_{S_i}^{(j+1)}$ and $ D_{S_{i+}}^{(j+1)}$ be the dominant transforms of $D_{S_i}^{(j)}$ and $D_{S_{i+}}^{(j)}$, whereas for $i=1,\dots k$ and $j=0,\dots k$ let $\Delta_{i+}^{(j+1)}$  be the strict transform of $\Delta_{i+}^{(j)}$ respectively. We prove the following:\\

$\mathit{Claim}$: For all $i$ and $j$ above, we have:\\
\begin{gather}\label{eq4}
\displaystyle\sum\limits _{i\in S_k=S} I_{\Delta_{i+}^{(j)}}= I_{D_{S_{j+}}^{(j)}}\cdot \displaystyle \prod \limits _{i=j+1}^{k} I_{D_{S_{i}}^{(j)}} + I_{\Delta_{i'+}^{(j)}}
\end{gather}
for any $i'\in S$.\\
\textit{Proof of claim}: The proof goes by induction on $j$. When $j=0$, we have:\\
\begin{gather*}
\displaystyle\sum\limits _{i\in S_k=S} I_{\Delta_{i+}^{(0)}}=\displaystyle\sum\limits _{i\in S_k=S} I_{\Delta_{i+}^{(Y_0)}}=\text{(by the definitions and Proposition \ref{prop1})}\\
\left(I_{D_{S_{0}}^{(0)}}\cap I_{D_{S_{1}}^{(0)}}\cap\dots \cap I_{D_{S_{k}}^{(0)}}\right)+ I_{\Delta_{i'+}^{(0)}}=\\
\displaystyle\prod\limits_{i=0}^{k} I_{D_{S_{i}}^{(0)}}+ I_{\Delta_{i'+}^{(0)}}
\end{gather*}
where the last equation follows from the fact that the $D_{S_{i}}^{(0)}$ intersect transversally (this is immediate since any collection of divisors $D_I$ in $\wmod$ intersects transversally by Theorem 2). The ideal $\displaystyle\prod\limits_{i=0}^{k} I_{D_{S_{i}}^{(0)}}+ I_{\Delta_{i'+}^{(0)}}$ doesn't change if we add the subideal $I_{\Delta_{i'+}^{(0)}}\cdot \displaystyle\prod\limits_{i=1}^{k} I_{D_{S_{i}}^{(0)}}$, so \\
\begin{gather*}
\displaystyle\prod\limits_{i=0}^{k} I_{D_{S_{i}}^{(0)}}+ I_{\Delta_{i'+}^{(0)}}=\\ 
\left(\displaystyle\prod\limits_{i=0}^{k} I_{D_{S_{i}}^{(0)}}+ I_{\Delta_{i'+}^{(0)}}\right)+I_{\Delta_{i'+}^{(0)}}\cdot \displaystyle\prod\limits_{i=1}^{k} I_{D_{S_{i}}^{(0)}}=\\
\left(I_{D_{S_{0}}^{(0)}}\cdot \displaystyle\prod\limits_{i=1}^{k} I_{D_{S_{i}}^{(0)}}+ I_{\Delta_{i'+}^{(0)}}\right)+ I_{\Delta_{i'+}^{(0)}}\cdot \displaystyle\prod\limits_{i=1}^{k} I_{D_{S_{i}}^{(0)}}=\\
\left(I_{D_{S_{0}}^{(0)}}+I_{\Delta_{i'+}^{(0)}}\right)\cdot \displaystyle\prod\limits_{i=1}^{k} I_{D_{S_{i}}^{(0)}}+I_{\Delta_{i'+}^{(0)}}=\text{(by the definitions)}\\
I_{D_{S_{0+}}^{(0)}} \cdot \displaystyle\prod\limits_{i=1}^{k} I_{D_{S_{i}}^{(0)}}+I_{\Delta_{i'+}^{(0)}}
 \end{gather*}
which proves the claim for $j=0$.\\
Now, assume the claim is true for $j$. We have\\
\begin{gather}
\displaystyle\sum\limits _{i\in S_k=S} I_{\Delta_{i+}^{(j+1)}}=(\text{by lemma \ref{lem2}.1(b))}\notag\\
\displaystyle\sum\limits _{i\in S_k=S}  I_{\pi_{j+1}^{-1}\Delta_{i+}^{(j)}}\cdot I_{\pi_{j+1}^{-1} D_{S_{j+}}^{(j)}}^{-1}=(\text{by (\ref{eq4})})\notag\\
\left(\displaystyle I_{\pi_{j+1}^{-1}D_{S_{j+}}^{(j)}}\cdot \prod \limits_{i=j+1}^{k} I_{\pi_{j+1}^{-1} D_{S_{i}}^{(j)}} +  I_{\pi_{j+1}^{-1}\Delta_{i'+}^{(j)}}\right)\cdot I_{\pi_{j+1}^{-1}D_{S_{j+}}^{(j)}}^{-1}
\end{gather}
Now, by lemma \ref{lem2}.1(b), we have
\begin{gather}
 I_{\pi_{j+1}^{-1}\Delta_{i'+}^{(j)}}=I_{\Delta_{i'+}^{(j+1)}}\cdot I_{\pi_{j+1}^{-1} D_{S_{j+}}^{(j)}}
\end{gather}
Also, by proposition \ref{prop2}, $D_{S_{i}}^{(j)},\,i\geq j+1$ is transversal to the center $D_{S_{j+}}^{(j)}$, hence, by lemma \ref{lem2}.1(a) \\
\begin{gather}
I_{\pi_{j+1}^{-1} D_{S_{i}}^{(j)}} = I_ {D_{S_{i}}^{(j+1)}} 
\end{gather}
Therefore, in view of (6) and (7), (5) reads \\
\begin{gather}
\left(\displaystyle I_{\pi_{j+1}^{-1}D_{S_{j+}}^{(j)}}\cdot \prod \limits_{i=j+1}^{k} I_{\pi_{j+1}^{-1} D_{S_{i}}^{(j)}} +  I_{\pi_{j+1}^{-1}\Delta_{i'+}^{(j)}}\right)\cdot I_{\pi_{j+1}^{-1}D_{S_{j+}}^{(j)}}^{-1}=\notag\\
\left(\displaystyle I_{\pi_{j+1}^{-1}D_{S_{j+}}^{(j)}}\cdot \prod \limits_{i=j+1}^{k} I_{ D_{S_{i}}^{(j+1)}} +  I_{\Delta_{i'+}^{(j+1)}}\cdot I_{\pi_{j+1}^{-1} D_{S_{j+}}^{(j)}}\right)\cdot I_{\pi_{j+1}^{-1}D_{S_{j+}}^{(j)}}^{-1}=\notag\\
\prod \limits_{i=j+1}^{k} I_{ D_{S_{i}}^{(j+1)}} +  I_{\Delta_{i'+}^{(j+1)}}=\notag\\
I_{ D_{S_{j+1}}^{(j+1)}}\cdot \prod \limits_{i=j+2}^{k} I_{ D_{S_{i}}^{(j+1)}} +  I_{\Delta_{i'+}^{(j+1)}}
\end{gather}
Adding the subideal $I_{\Delta_{i'+}^{(j+1)}}\cdot \prod \limits_{i=j+2}^{k} I_{ D_{S_{i}}^{(j+1)}}$ to (8) doesn't change anything so the ideal of (8) is equal to 
\begin{gather}
\left(I_{ D_{S_{j+1}}^{(j+1)}}\cdot \prod \limits_{i=j+2}^{k} I_{ D_{S_{i}}^{(j+1)}} +  I_{\Delta_{i'+}^{(j+1)}}\right)+I_{\Delta_{i'+}^{(j+1)}}\cdot \prod \limits_{i=j+2}^{k} I_{ D_{S_{i}}^{(j+1)}}=\notag\\
\left(I_{ D_{S_{j+1}}^{(j+1)}}+I_{\Delta_{i'+}^{(j+1)}}\right)  \cdot \prod \limits_{i=j+2}^{k} I_{ D_{S_{i}}^{(j+1)}} + I_{\Delta_{i'+}^{(j+1)}}=\text{(by Proposition \ref{prop3})}\notag\\
I_{ D_{S_{j+1+}}^{(j+1)}}\cdot \prod \limits_{i=j+2}^{k} I_{ D_{S_{i}}^{(j+1)}} + I_{\Delta_{i'+}^{(j+1)}}\notag
 \end{gather}
 which proves the claim for $j+1$, hence the proof of the  claim is now complete.\\
 Setting $j=k$ in the claim shows (since  $\prod \limits _{i=j+1}^{k} I_{D_{S_{i}}^{(j)}}$ is vacuous) that \\
 \begin{gather*}
 \displaystyle\sum\limits _{i\in S_k=S} I_{\Delta_{i+}^{(k)}}=I_{D_{S_{k+}}^{(k)}}+ I_{\Delta_{i'+}^{(k)}}= I_{D_{S_{k+}}^{(k)}}
   \end{gather*}
   where in the last equality we use the fact that $\Delta_{i'+}^{(k)}\supset D_{S_{k+}}^{(k)}$, since $i'\in S_k=S$. By lemma \ref{lem2}(2), we see that the strict transforms of the $\Delta_{i+}^{(k)}$ do not intersect simultaneously in $X_{k+1}$, that is, the iterated strict transforms of $\Delta_{i+}^{(0)} =\Delta_{i+}^{(Y_0)},\, i\in S$ do not intersect simultaneously in the open cover $\{U_C\}_{C\in\cc}$ of stage $(k,S)$, hence 
   \begin{center}
   $\displaystyle \bigcap\limits _{i\in S, |S|=n-k}\Delta_{i+}^{(S+)}=\emptyset$ 
   \end{center}
  in stage $(k,S)$. This shows part (i) of the Proposition.\\
  
  Part (ii) now follows from (i). Let $S,T\subset N$ such that $\Delta_S, \Delta_T \in \wbuild$ and $|S|=|T|=n-k-1, S\neq T$.  By definition, we have $\Delta_{i+}^{(Y_0)}\supset D_{S+}^{(Y_0)}$ for all $i\in S$ and $\Delta_{i+}^{(Y_0)}\supset D_{T+}^{(Y_0)}$ for all $i\in T$. Therefore we deduce the inclusions of their iterated strict transforms $\Delta_{i+}^{(Y_{k+1})}\supset D_{S+}^{(Y_{k+1})}$ for all $i\in S$ and $\Delta_{i+}^{(Y_{k+1})}\supset D_{T+}^{(Y_{k+1})}$  for all $i\in T$ respectively. This, in turn, implies that \\
  \begin{gather}
D_{S+}^{(Y_{k+1})}\cap D_{T+}^{(Y_{k+1})} \subset \displaystyle \bigcap\limits _{i\in S\cup T}\Delta_{i+}^{(Y_{k+1})} ,
\end{gather}    
   
     \vspace{0.1in}
However, the assumptions on $S$ and $T$ imply that $|S\cup T|\geq n-k$. Let $U\subset S\cup T$ be such that $|U|=n-k$ and $\Delta_U \in \wbuild$. By part (i), we have $\displaystyle \bigcap\limits _{i\in U}\Delta_{i+}^{(U+)}=\emptyset$. Since $Y_{k+1}$ is obtained from stage $(k,U)$ by a finite sequence of blowups we have that the strict transforms of the $\Delta_{i+}^{(U+)}$ in $Y_{k+1}$ also have empty intersection, that is:\\
\begin{gather}
\displaystyle \bigcap\limits _{i\in U}\Delta_{i+}^{(Y_{k+1})}=\emptyset 
\end{gather}
Finally, (9) and (10) together imply that\\
\begin{gather*}
D_{S+}^{(Y_{k+1})}\cap D_{T+}^{(Y_{k+1})} \subset \displaystyle \bigcap\limits _{i\in S\cup T}\Delta_{i+}^{(Y_{k+1})} \subset \displaystyle \bigcap\limits _{i\in U}\Delta_{i+}^{(Y_{k+1})}=\emptyset 
\end{gather*}
that is $D_{S+}^{(Y_{k+1})}\cap D_{T+}^{(Y_{k+1})}=\emptyset$, hence the proof of part (ii).\\

The proof of part (iii) is identical to the proof of (ii) and is therefore omitted. \qed\\

\textit{Proof of Theorem 3}: By the above Proposition, it follows that the fibers of $\wmod^+\rightarrow \wmod$ are the $\mathcal{A}$ stable degenerations described in section 2.2.1 with $n$ sections induced by the strict transforms of the $\Delta_{i+}^{(Y_0)}\subset \wmod\times X$ under the sequence of blowups $\wmod^+\rightarrow \wmod\times X$ and vice versa. Since the fibers are then equidimensional and the varieties $\wmod^+$ and $\wmod$ are nonsingular, we conclude that  $\wmod^+\rightarrow \wmod$ is flat. \qed



 \end{subsection}
 \vspace{0.1in}
 
 \begin{subsection}{Functorial Description}\label{sec2.4}Let $\wbuild$ as above and $S\subset N$ such that $\Delta_S\in \wbuild$. 
 Now let $F:(Sch/k)\rightarrow (Set)$ be the functor that  to every scheme $H$ assigns the set of pairs of morphisms (of schemes and sheaves) \\
 \begin{center}
$ \{f:H\rightarrow X^n \, ,f^*I_{\Delta_S} \twoheadrightarrow {\mathcal{L}}_S\}$
 \end{center}
 \bigskip
for all $S\subset N$ such that $\Delta_S\in \wbuild$ , with ${\mathcal{L}}_S$  invertible and the compatibility condition that if $T\subset S\subset N$ there exists a (necessarily unique) morphism ${\mathcal{L}}_T\rightarrow   {\mathcal{L}}_S$ such that the diagram \\
 
\centerline{
 \xymatrix{
f^*I_{\Delta_T} \ar[d] \ar@{->>}[r] & {\mathcal{L}}_T \ar[d] \\
f^*I_{\Delta_S} \ar@{->>}[r] & {\mathcal{L}}_S
} }
commutes.
\begin{thm}\label{thm4}(cf. \cite[Theorem 4] {FM}) The functor $F$ is representable by the variety $\wmod$. 
\end{thm}
\textbf{Proof}: First, let us assume that there exists a morphism $g: H \rightarrow \wmod$. Composing with the blowup morphism $\pi: \wmod \rightarrow X^n$, we get a morphism $f=\pi \circ g: H \rightarrow X^n$.\\

\centerline{
\xymatrix{
&\wmod \ar[d] ^\pi \\
H \ar [ru]^g \ar[r]_f & X^n}
}
\bigskip
Now, set ${\mathcal{L}}_S:=g^* I_{\pi^{-1}\Delta_S}$. By theorem \ref{thm1}(3), $\wmod$ is isomorphic to the iterated blowup along the total inverse images of the ideal sheaves of the diagonals that belong to $\wbuild$. Therefore, $I_{\pi^{-1}\Delta_S}$ is invertible, which implies that ${\mathcal{L}}_S$ is invertible as well. Moreover, we have a canonical surjection \\
\begin{center}
$\pi^*I_{\Delta_S} \twoheadrightarrow I_{\pi^{-1}\Delta_S}$
\end{center}
which further induces a surjection \\
\begin{center}
$f^*I_{\Delta_S}\iso g^*\pi^*I_{\Delta_S} \twoheadrightarrow g^*I_{\pi^{-1}\Delta_S}={\mathcal{L}}_S$.
\end{center}

\bigskip

 If $T\subset S$ then $\Delta_T \in \wbuild$, so $I_{\pi^{-1}\Delta_T}$ is also invertible. We  obtain the following commutative diagram of sheaves on $H$:\\
 
 \centerline{
 \xymatrix{
 \pi^*I_{\Delta_T} \ar[d] \ar@{->>}[r] & I_{\pi^{-1}\Delta_T}\ar[d]\\
  \pi^*I_{\Delta_S} \ar@{->>}[r]  & I_{\pi^{-1}\Delta_T}
  }}  
  \bigskip
  We then apply $g^*$ to the above diagram to get the desired commutative diagram :\\
  
   \centerline{
 \xymatrix{
f^*I_{\Delta_T} \iso g^*\pi^*I_{\Delta_T} \ar[d] \ar@{->>}[r] & g^*I_{\pi^{-1}\Delta_T}={\mathcal{L}}_T\ar[d]\\
f^*I_{\Delta_S} \iso g^* \pi^*I_{\Delta_S} \ar@{->>}[r]  & g^*I_{\pi^{-1}\Delta_S}= {\mathcal{L}}_S  }}  
 
  \bigskip
  Conversely, let us assume that there exists a morphism $f:H\rightarrow X^n$ and compatible invertible quotients $f^*I_{\Delta_S} \twoheadrightarrow {\mathcal{L}}_S$ for all $S\subset N$ such that $\Delta_S \in \wbuild$ .  We will show that we can lift $f$ to $\wmod$ by induction. First, observe that $\wmod$ may be obtained as a sequence of blowups along iterated dominant transforms of the elements of $\wbuild$, this time written in a different order as follows:\\
 
  Consider the set  \\
  \begin{center}
  ${\mathcal{G}}^{FM}:\{\Delta_{12}; \Delta_{123}; \Delta_{13}, \Delta_{23}; \Delta_{1234}; \Delta_{124}, \Delta_{134}, \Delta_{234}; \Delta_{14},  \Delta_{24},   \Delta_{34};\dots \} $
  \end{center}
  \bigskip
  This is the set of blowup centers used in the Fulton and Macpherson compactification \cite{FM} (also, see example 1). The strict transforms of the diagonals separated by a comma are disjoint, so for these particular centers the order that the blowups are taken does not matter. Now, consider the set obtained from ${\mathcal{G}}^{FM}$ by deleting the diagonals indexed by coordinates whose sum of weights is less than or equal to 1:\\
  \begin{center}
  $\mathcal{G}_{\mathcal{A}}^{FM}:=\mathcal{G}\setminus \{ \Delta_S | S\subset N \, \text{and} \sum \limits _{i\in S} a_i \leq 1\}$
  \end{center}
  \bigskip
Now write its elements preserving the above order of ${\mathcal{G}}^{FM}$. We make the following claim: \\

\textit{Claim}: The set $\mathcal{G}_{\mathcal{A}}^{FM}$ satisfies the assumption of Theorem \ref{thm1}(2). \\
\textit{Proof of Claim}: Consider the set of the first $i$ elements $\mathcal{G}_{\mathcal{A},i}^{FM}$ and pick an arbitrary subset of the latter consisting of $k$ elements $\{\Delta_{I_1}, \Delta_{I_2},\dots , \Delta_{I_k}\}$, where $I_1<I_2\dots <I_k$ in the above order. Now, take the intersection $S_k:=\Delta_{I_1}\cap \Delta_{I_2}\cap\dots \cap \Delta_{I_k}$. We need to show that the minimal elements of $\mathcal{G}_{\mathcal{A},i}^{FM}$ containing $S_k$ intersect transversally along $S_k$. Let $l$ be the greatest integer contained in $I_k$. We partition the above set of indices in two subsets  $\{\Delta_{I_1}, \Delta_{I_2},\dots , \Delta_{I_{k'}}\}$ and $\{\Delta_{I_{k'+1}}, \Delta_{I_{k'+2}},\dots , \Delta_{I_k}\}$ such that the greatest integer contained in each of the $I_{k'+1}, I_{k'+2},\dots I_k$ is $l$, whereas the greatest integer contained in $I_{k'}$ is $l-1$. The first set could be empty; in this case the claim follows trivially, so we assume that it is non empty. Then we have $\Delta_{I_{k'+1}}\cap \Delta_{I_{k'+2}}\cap\dots \cap \Delta_{I_k}=\Delta_{I\cup\{l\}}$, where $I\subset\{1,2,\dots l-1\}$. Now consider the intersection $\Delta_{I_1}\cap \Delta_{I_2}\cap\dots \cap \Delta_{I_{k'}}\cap \Delta_I$. As in the proof of Lemma \ref{lem1}, the latter can be written uniquely as $\Delta_{J_1}\cap\Delta_{J_2}\dots\cap \Delta_{J_{k''}}$, where $k''\leq k'$ and for  $j=1,\dots k''$, the sets $J_j$ are unions of some of the $I_m$ and $I$ ($m=1,\dots k'$) and also pairwise disjoint.
Then, there exists a unique set, say $J_s$, among the $J_j$ above that contains $I$. We have\\

\begin{gather*}
S_k:=(\Delta_{I_1}\cap \Delta_{I_2}\cap\dots \Delta_{I_{k'}})\cap(\Delta_{I_{k'+1}}\dots \cap \Delta_{I_k})\\
=(\Delta_{J_1}\cap\Delta_{J_2}\dots\cap \Delta_{J_{k''}}) \cap \Delta_{I\cup\{l\}}\\
=(\Delta_{J_1}\cap\Delta_{J_2}\dots\Delta_{J_{s-1}}\cap\Delta_{J_{s+1}}\cap \dots \Delta_{J_{k''}}) \cap (\Delta_{J_s}\cap\Delta_{I\cup\{l\}})\\
=\Delta_{J_1}\cap\Delta_{J_2}\dots\Delta_{J_{s-1}}\cap\Delta_{J_{s+1}}\cap \dots \Delta_{J_{k''}} \cap \Delta_{J_s\cup\{l\}}\\
\end{gather*}

\bigskip

Since the $J_j$ are pairwise disjoint subsets of $\{1,2,\dots, l-1\}$ we see that the indices appearing in the last equation above are also pairwise disjoint. Moreover, each of the $J_j$ with $j\neq s$ contains some $I_m$ where $m\leq k'$, so \\
\begin{center}
 $\sum \limits _{i\in J_j} a_i \geq \sum \limits _{i\in I_m} a_i>1$
\end{center}
\vspace{0.1in}
In addition, since the $J_j$ above with $j\neq s$ are unions of the $I_m\subset \{1,2,\dots, l-1\}$ the maximum integer contained in each of them is less than or equal to $l-1$.  This implies that $\Delta_{J_j}<\Delta_{I_{k'+1}}$ for all $j\neq s$, which in turn implies $\Delta_{J_j}\in \mathcal{G}_{\mathcal{A},i}^{FM}$ for all $j\neq s$.
     
 \vspace{0.1 in}
 
Similarly we see that $\Delta_{J_s\cup\{l\}}\in  \mathcal{G}_{\mathcal{A},i}^{FM}$. Therefore, the minimal elements in $\mathcal{G}_{\mathcal{A},i}^{FM}$ that contain $S_k$ are precisely the sets $ \Delta_{J_1}, \Delta_{J_2}, \Delta_{J_{s-1}}, \Delta_{J_{s+1}}, \dots \Delta_{J_{k''}},\Delta_{J_s\cup\{l\}}$, which are seen to intersect transversally -as varieties- since the indices they are indexed by are disjoint. \\

In view of the claim, we deduce, by Theorem \ref{thm1}(2), that $\wmod$ is isomorphic to the blowup along the iterated dominant transforms of the elements of $\mathcal{G}_{\mathcal{A}}^{FM}$ in the order described above. The proof is then essentially identical to the proof in \cite[theorem 4]{FM}.
\qed
 \end{subsection}
 \begin{subsection}{Reduction and Forgetful morphisms}\label{sec2.5}
 The results of this section are the analogues of \cite[Section 4]{Hassett}. See also \cite{AG}, \cite{BayerManin}.

 \begin{thm}\label{thm5}\label{thm5} Let $\mathcal{A}:=\{a_1,a_2, \dots,a_n\} $ and  $\mathcal{B}:=\{b_1,b_2, \dots,b_n\} $ be two sets of rational numbers such that $0<b_i\leq a_i\leq1$ for all $i=1,2,\dots n$. There exists a natural blowup (reduction) morphism \\
 \begin{center}
$\rho_{\mathcal{B},\mathcal{A}}: \wmod \rightarrow \wmodb$
 \end{center}
 \vspace{0.1in}
 At the level of $k$ points, the above morphism reassigns the weights of the sections of an  $\mathcal{A}$ stable degeneration and then successively collapses all components that are unstable with respect to $\mathcal{B}$. 
 \end{thm}
 \textbf{Proof}: As in section 2.2, let \\
\begin{center}
$\mathcal{G}:=\{ \Delta_{12\dots n}; \Delta_{12\dots (n-1)},\dots ,\Delta_{23\dots n}; \dots;   \Delta_{12},\dots   ,\Delta_{(n-1)n} \} $
\end{center}
\vspace{0.1in}

and
\begin{center}
$\mathcal{G}_{\mathcal{B}}:=\mathcal{G}\setminus \{ \Delta_S | S\subset N \, \text{and}\, \sum \limits _{i\in S} b_i \leq 1\}$

\end{center}
and write the elements of the latter in ascending dimension order.
  By Lemma \ref{lem1}, $\wbuildb$ is a building set. By the hypothesis, $\wbuildb\subset\wbuild$. For ease of notation, denote the ideal sheaves of the $\Delta_S \in \wbuild$ by $I_1, I_2,\dots I_k$ in order preserving bijection with the $\Delta_S$ in $\wbuild$. Similarly, let $\{I_{i_1}, I_{i_2}, \dots I_{i_l}\}\subset \{I_1, I_2, \dots I_k\}$ be the set of ideal sheaves of the $\Delta_S \in \wbuildb$, listed again in order preserving bijection with the $\Delta_S$ in $\wbuildb$. We have, \\
  \begin{gather*}
 \wmod\iso  Bl_{I_k}\dots Bl_{I_2}Bl_{I_1} X^n \, (\text{theorem \ref{thm1}(3))}\\
 \iso Bl_{I_k}\dots \left (Bl_{I_{i_l}}\dots Bl_{I_{i_2}}Bl_{I_{i_1}}X^n \right) \, (\text{\cite[Lemma 3.2]{Li}}) \\
 \iso Bl_{I_k} \dots \wmodb 
  \end{gather*} 
  therefore there is a natural blowup map along centers corresponding to $\Delta_S \in \wbuild \setminus \wbuildb$\\
  \begin{center}
  $\rho_{\mathcal{B},\mathcal{A}}: \wmod \rightarrow \wmodb$
  \end{center}
  \qed
  \vspace{0.1in}
  
  With notation as in Theorem 5, we have:
  
   \begin{prop}\label{prop5} Let $\mathcal{A}:=\{a_1,a_2, \dots,a_n\} $,  $\mathcal{B}:=\{b_1,b_2, \dots,b_n\} $ and  $\mathcal{C}:=\{c_1,c_2, \dots,c_n\} $ be sets of rational numbers such that $0<c_i\leq b_i\leq a_i\leq1$ for all $i=1,2,\dots n$. Then: \\
 \begin{center}
$\rho_{\mathcal{C},\mathcal{A}}= \rho_{\mathcal{C},\mathcal{B}}\circ \rho_{\mathcal{B},\mathcal{A}}$
 \end{center}  
 \end{prop}
 \textbf{Proof}: With notation as in the proof of Theorem 5 above, we have $\wbuildc\subset \wbuildb \subset \wbuild$. We also see that $\rho_{\mathcal{B},\mathcal{A}}:\wmod\rightarrow \wmodb$ is the blowup along centers corresponding to $\wbuild \setminus \wbuildb$, whereas  $\rho_{\mathcal{C},\mathcal{B}}:\wmodb\rightarrow\wmodc$ is the blowup along centers  corresponding to  $\wbuildb \setminus \wbuildc$, hence their composition is the blowup along centers corresponding to $\wbuild \setminus \wbuildc$. This is precisley the morphism $\rho_{\mathcal{C},\mathcal{A}}:\wmod\rightarrow \wmodc$, as we saw in the proof of Theorem 5. \qed

\bigskip
 
  \begin{thm}\label{thm6}
  Let $\mathcal{A'}:=\{a_{i_1},a_{i_2},\dots a_{i_r}\}$ be a subset of $\mathcal{A}$ such that $r\leq n$. There exists a natural forgetful morphism\\
 \begin{center}
$\phi_{\mathcal{A},\mathcal{A'}}: \wmod \rightarrow \wmods $
\end{center}
\vspace{0.1in}
 At the level of $k$ points, the above morphism successively collapses all components of an $\mathcal{A}$ stable degeneration that are unstable with respect to $\mathcal{A'}$. \\
    \end{thm}
\textbf{Proof}: We begin with the following general \textit{remark}: let $\{I_1, I_2,\dots, I_m\}$ be a collection of subsets of the set $\{i_1, i_2, \dots i_r\}$. Then we have\\
\begin{gather}
\left(Bl_{\Delta_{I_m}}\dots Bl_{\Delta_{I_1}} X^r\right) \times_{X^r} X^n \iso Bl_{\Delta_{I_m}}\dots Bl_{\Delta_{I_1}} X^n 
\end{gather} 
\vspace{0.1in}

We explain the notation: on the left hand side of the above isomorphism each diagonal ${\Delta_{I_j}}$ is considered a subvariety inside $X^r$, whereas on the right hand side we consider the corresponding diagonal in $X^n$, that is we consider $I_j$ as a subset of $N$. On both sides, $X^r$ and $X^n$ respectively are successively blown up along the \textbf{total inverse images} of the $\Delta_{I_k}$. Now, the isomorphism (11) is a direct consequence of the fact that the formation of the blowup commutes with smooth base change. \\
 
Next, consider $\{i_1, i_2, \dots i_r\}$ as a subset of $N$. Let \\
\begin{center}
$\mathcal{G}:=\{ \Delta_{12\dots n}; \Delta_{12\dots (n-1)},\dots ,\Delta_{23\dots n}; \dots;   \Delta_{12},\dots   ,\Delta_{(n-1)n} \} $
\end{center}
\vspace{0.1in}

and
\begin{center}
$\mathcal{G}_{\mathcal{A'}}:=\mathcal{G}\setminus \{ \Delta_S | S\subset N \, \text{and} \sum \limits _{i_k\in S} a_{i_k} \leq 1\}$

\end{center}
and write the elements of the latter in ascending dimension order.
 As in the proof of Theorem \ref{thm5}, $\wbuilds$ is a building set (of its induced arrangement of subvarieties in $X^n$) such that $\wbuilds \subset \wbuild$ and if $\{J_1, J_2, \dots J_l\}$ 
 are the ideal sheaves of the $\Delta_S \in \wbuilds $, then there is a natural blowup map along centers corresponding to $\Delta_S \in \wbuild \setminus \wbuilds$\\
 \begin{gather}
 \wmod \rightarrow Bl_{J_l}\dots Bl_{J_1} X^n
 \end{gather}

\bigskip 
 Now, we go back to the remark of the first paragraph. With notation as in the first paragraph, set $m=l$ and take $\Delta_{I_{k}}$ (where $\Delta_{I_k}$ is considered as a subvariety of $X^n$) to be the diagonal whose ideal sheaf is $J_k,\, k=1,\dots l$. We then have\\
 \begin{gather}
 Bl_{J_l}\dots Bl_{J_1} X^n \iso  \left(Bl_{J_l}\dots Bl_{J_1} X^r\right) \times_{X^r} X^n
  \end{gather}
 
  \bigskip
 where the $J_k$ in the first isomorphism of (13) are the ideal sheaves of $\Delta_{I_k}$  considered as  subvarieties of $X^r$. 
 
 \bigskip 
 By  Theorem \ref{thm1}(3), we have $Bl_{J_l}\dots Bl_{J_1} X^r\iso \wmods$, hence
  \begin{gather}
  \left(Bl_{J_l}\dots Bl_{J_1} X^r\right) \times_{X^r} X^n \iso \wmods \times_{X^r} X^n
\end{gather}   
  
  \bigskip
  
 Therefore composing the morphisms in (12), (13) and (14) we get a morphism  $\wmod \rightarrow \wmods \times_{X^r} X^n$. Composing further with the projection $\wmods \times_{X^r} X^n \rightarrow \wmods$ we obtain the required morphism  \\
 \begin{center}
     $\phi_{\mathcal{A},\mathcal{A'}}: \wmod \rightarrow \wmods $
 \end{center}
     \qed
   
       \end{subsection}
  
  
  
\end{section}
 \begin{section}{Relative stable degenerations}\label{sec3}
 \begin{subsection}{The relative weighted compactification}\label{sec3.1}
 Let $f: X\rightarrow Y $ be a smooth morphism between two smooth varieties  $X$ and $Y$. Moreover, let $\mathcal{A}:=\{a_1,a_2, \dots,a_n\} $ be a set of rational numbers such that\\
\renewcommand{\labelitemi}{$\bullet$}
\begin{itemize}
\item $0<a_i\leq 1 , i=1,2,\dots,n$ 
\end {itemize}
\vspace{0.1in}
Let

\begin{center}
$\mathcal{G}^{rel}:=\{ \Delta_{12\dots n}; \Delta_{12\dots (n-1)},\dots ,\Delta_{23\dots n}; \dots;  \Delta_{12},\dots   ,\Delta_{(n-1)n} \} $
\end{center}
\bigskip
be the set of diagonals of the $n$-fold fibered product $X_Y^n:=\underbrace{X\times_Y\times\dots \times_Y X}_{n factors}$ listed in ascending dimension order.
\vspace{0.1in}

Now let 
\begin{center}
$\wbuildrel:=\{ \Delta_S\subset X_Y^n | S\subset N \, \text{and}\, \sum \limits _{i\in S} a_i > 1\}$

\end{center}
 and list its elements in ascending dimension order. \\
 
 \vspace{0.1in}
\begin{mydef} The relative weighted configuration space with respect to $\mathcal{A}$ is the complement of the union of the diagonals that belong to $\wbuildrel$ in $X_Y^n$\\
 \begin{center}
 $F_\mathcal{A}^{rel}(X/Y,n):= X_Y^n\setminus \bigcup\limits_ {\Delta_I \in \wbuildrel} \Delta_I$ . \\
  \end{center} 
\end{mydef}
 The relative weighted configuration space may be viewed as the parameter space of $n$ labeled points on the fibers of $f$ carrying weights $a_i$ subject to the following condition:\\
\renewcommand{\labelitemi}{$\bullet$}
\begin{itemize}
\item for any set of labels $S\subset N$of coincident points we have $\sum \limits _{i\in S} a_i \leq 1\ $.
\end {itemize}

 We construct a weighted compactification of the above configuration space as the space of $\mathcal{A}$ stable degenerations (see Section 2) on the fibers of $f$. \\
  \begin{mydef}The relative weighted compactification of the relative configuration space $F_\mathcal{A}^{rel}(X/Y,n)$ is  the closure of the image of the diagonal embedding\\
 
  \begin{center}
 
$ X_Y^n\setminus \bigcup\limits_ {\Delta_I \in \wbuildrel} \Delta_I\xhookrightarrow\ X_Y^n \times\displaystyle \prod \limits_ {\Delta_I\in\ \wbuildrel} Bl_{\Delta_I} X_Y^n$

\end {center}
  We denote the relative weighted compactification by $\wmodrel$.
\end{mydef}
\vspace{0.1in}
\textbf{Remark}:In case $\sum_{i=1}^n a_i\leq1$, then $\wbuildrel$ is empty and $\wmodrel$ is trivially equal to $X_Y^n$.\\

\textbf{Remark}: Since the relative configuration space of $n$ \textit{distinct} labeled points $X_Y^n\setminus \bigcup\limits_ { |I| \geq 2} \Delta_I $ is an open subset of $ X_Y^n\setminus \bigcup\limits_ {\Delta_I \in \wbuildrel} \Delta_I$,  the closure of 
\begin{center}
 $X_Y^n\setminus \bigcup\limits_ { |I| \geq 2} \Delta_I\xhookrightarrow\ X_Y^n\setminus \bigcup\limits_ {\Delta_I \in \wbuildrel} \Delta_I\xhookrightarrow\ X_Y^n \times\displaystyle \prod \limits_ {\Delta_I\in\ \wbuildrel} Bl_{\Delta_I} X_Y^n$ 
 \end{center}
 is the same as  the closure of\\
  \begin{center}
 
$ X_Y^n\setminus \bigcup\limits_ {\Delta_I \in \wbuildrel} \Delta_I\xhookrightarrow\ X_Y^n \times\displaystyle \prod \limits_ {\Delta_I\in\ \wbuildrel} Bl_{\Delta_I} X_Y^n$.

\end {center}
Therefore, the variety $\wmodrel$ may be viewed as a compactification of the relative configuration space $F^{rel}(X/Y,n)$ of $n$ \textit{distinct} labeled points on the fibers of $X\rightarrow Y$ and consequently as a `weighted' generalization of the Fulton Macpherson compactification to the relative context.

 \begin{prop}\label{prop6}The relative weighted compactification $X_Y^{\mathcal{A}}[n]$ is isomorphic to the iterated blowup of the $n$-fold fibered product $X_Y^n$ along the dominant transforms of the diagonals in $\wbuildrel$ in ascending dimension order.
 
  \end{prop}
  \textbf{Proof}: The statement is Zariski local, so we may assume $Y$ is irreducible after base changing to the (open) irreducible components of $Y$(note that $Y$ is smooth). Moreover, it suffices to prove the statement after restricting to the irreducible components of $X_Y^n$; those are open and disjoint, since $X_Y^n$ is a smooth variety. So, let $U$ be an (open) irreducible component of $X_Y^n$. For any $I\subset N$ such that $|I|\geq2$, define $\Delta_I^U:=\Delta_I\cap U$ and set\\
  \begin{center}
  ${\wbuildrelU}:=\{\Delta_I\cap U |\Delta_I \in\wbuildrel \, \text{and}\, \Delta_I^U\neq\emptyset\}$
  \end{center}
  \bigskip
  and, as usual, write the elements of the above set in ascending dimension order.\\
  
  \textit{Claim}: $\wbuildrelU$ is a building set of subvarieties of $U$.\\
  \textit{Proof of Claim}: First, observe that the schemes $\Delta_I^U \in\wbuildrelU$ are irreducible subvarieties of $U$. Indeed, each diagonal $\Delta_I$ may be considered as an embedding $X_Y^{I'}\xhookrightarrow\  X_Y^N$, where $I'=(N\setminus I)\cup\{ i\}$ for any element $i\in I$, induced by the identity morphism on $X_Y^{N\setminus I}$ and the diagonal embedding $X^{\{i\}}\xhookrightarrow\ X_Y^I$ over $Y$.  This embedding is then  a section $X_Y^{I'}\xhookrightarrow\  X_Y^N\xrightarrow{p_{I'}} X_Y^{I'}$ of the projection $p_{I'}$ onto the coordinates of $X_Y^N$ labeled by $I'$ . Both $X_Y^{I'}$ and $X_Y^N$ are disjoint unions of their irreducible(connected) components. If $\Delta_I\cap U \neq \emptyset$, then any component of $\Delta_I=X_Y^{I'}$ either is contained in $U$ or is disjoint from it.  Moreover we cannot have two (disjoint) components of $X_Y^{I'}$ contained in $U$; this would imply that the image of the irreducible component $U$ under $p_I$  meets those two disjoint components, which is absurd. Therefore, if $\Delta_I\cap U\neq \emptyset$, then $\Delta_I^U$ is equal to the unique component of $\Delta_I=X_Y^{I'}$ contained in $U$, hence it is irreducible.\\
  
  To conclude the proof of the claim consider a nonempty intersection of elements of $\wbuildrelU$, say $S_k=\Delta_{I_1}^U\cap \Delta_{I_2}^U\cap\dots \cap \Delta_{I_k}^U$, where $\Delta_I\in \wbuildrel$. We need to show that the minimal elements of $\wbuildrelU$ that contain the above intersection intersect transversally and that their intersection is $S_k$. By the proof of Lemma \ref{lem1}, we  can write the  intersection $\Delta_{I_1}\cap \Delta_{I_2}\cap\dots \cap \Delta_{I_k}$ as an intersection of diagonals $\Delta_{J_1}\cap \Delta_{J_2}\cap\dots \cap \Delta_{J_k'}$ in $X_Y^n$, where $k'\leq k$ and the $J_j$ are unions of some of the $I_i$ and pairwise disjoint. Since $S_k$ is nonempty, restricting to $U$ gives the relation $S_k=\Delta_{I_1}^U\cap \Delta_{I_2}^U\cap\dots \cap \Delta_{I_k}^U=\Delta_{J_1}^U\cap \Delta_{J_2}^U\cap\dots \cap \Delta_{J_{k'}}^U$ in $U$, where the $\Delta_{J_j}^U$ are nonemtpy for $j=1,\dots,k'$. The sum of the weights of the coordinates of $\Delta_ {J_j}$ is greater than 1, because each $J_j$ contains some $I_i$. Also, the $\Delta_{J_j}^U$ intersect transversally along $S_k$, since they are indexed by disjoint sets. What remains to be verified in this case is that the $\Delta_{J_j}^U$ are the minimal elements of $\wbuildrelU$ that contain $S_k$. To see this, let $\Delta _S\in \wbuildrel$ be such that\\
  \begin{gather}\label{eq15}
 \Delta_S^U\supseteq \Delta_{J_1}^U\cap \Delta_{J_2}^U\cap\dots \cap \Delta_{J_{k'}}^U
 \end{gather}
 \bigskip
 Observe that $S$ must be contained in the union $\mathcal{I}:=\bigcup\limits_{j=1}^{k} J_j$. For assume the contrary; we then distinguish two cases: \\
 \begin{enumerate}[(a)]
 \item{ $S\cap\mathcal{I}\neq \emptyset$.} Then S can be written as a union of two nonempty disjoint subsets $S=S'\cup S''$, where $S'\subset\mathcal {I}$ and $S''\cap \mathcal{I}=\emptyset$. Assume, without loss of generality, that $J_1$ intersects $S'$ and set $S_0:=(S'\cap J_1)\cup S''$. Then, $\Delta_{S_0}^U\supseteq \Delta_S^U$, so by (15) we get\\
 \begin{gather*}
  \Delta_{S_0}^U\cap \Delta_{J_1}^U\cap \Delta_{J_2}^U\cap\dots \cap \Delta_{J_{k'}}^U=\Delta_{J_1}^U\cap \Delta_{J_2}^U\cap\dots \cap \Delta_{J_{k'}}^U
  \end{gather*}
  which, since $\Delta_{S_0}^U\cap \Delta_{J_1}^U=\Delta_{J_1\cup S''}^U$, implies that\\
\begin{gather*}
   \Delta_{J_1\cup S''}^U\cap \Delta_{J_2}^U\cap\dots \cap \Delta_{J_{k'}}^U=\Delta_{J_1}^U\cap \Delta_{J_2}^U\cap\dots \cap \Delta_{J_{k'}}^U
  \end{gather*}
 
 \vspace{0.1in} 
 
 All indices appearing on both sides of the above equation are pairwise disjoint so the intersections on both sides are transversal. Therefore the codimension of both sides in $U$ is equal to the sum of the codimensions of each of the subvarieties appearing in the intersections on each side. However, we saw above that the (nonempty) intersection of $U$ with any diagonal is the unique component of the diagonal contained in $U$. Since the diagonals are pure dimensional, we have  $dim \Delta_{J_1\cup S''}^U > dim \Delta_{J_1}^U$, so we get a contradiction.
\vspace{0.1in}
  \item{$S\cap\mathcal{I}= \emptyset$}.  In this case, again by equation (\ref{eq15}), we obtain the equality:
  \begin{gather*}
    \Delta_S^U\cap \Delta_{J_1}^U\cap \Delta_{J_2}^U\cap\dots \cap \Delta_{J_{k'}}^U=\Delta_{J_1}^U\cap \Delta_{J_2}^U\cap\dots \cap \Delta_{J_{k'}}^U
   \end{gather*}
  
   \vspace{0.1in}
  
   Since all indices appearing on each side are disjoint, each of the intersections appearing in the above equation is transversal. This implies, in particular, that $cod(\Delta_S^U,U)=0$, which is a contradiciton.
     \end{enumerate}
 
 \vspace{0.1in}
 Next, we show that there must be a unique set among the $J_j\, ,j=1,\dots k'$ containing $S$. Indeed, assume that S intersects at least two sets of indices among the $J_j$, and suppose without loss of generality that those are the sets $J_1$ and $J_2$. Let $S^*=S\cap (J_1\cup J_2)$. Then $\Delta_{S^*}^U\supseteq \Delta_S^U$, so, by equation (\ref{eq15}), we conclude that \\
 \begin{gather*}
  \Delta_{S^*}^U\cap \Delta_{J_1}^U\cap \Delta_{J_2}^U\cap\dots \cap \Delta_{J_{k'}}^U=\Delta_{J_1}^U\cap \Delta_{J_2}^U\cap\dots \cap \Delta_{J_{k'}}^U
  \end{gather*}

  \vspace{0.1in}
  Now, by the definition of $S^*$, we have $\Delta_{S^*}^U\cap \Delta_{J_1}^U\cap \Delta_{J_2}^U=\Delta_{J_1\cup J_2}^U$, hence 
           \begin{gather*}
  \Delta_{J_1\cup J_2}^U\cap\dots \cap \Delta_{J_{k'}}^U=\Delta_{J_1}^U\cap \Delta_{J_2}^U\cap\dots \cap \Delta_{J_{k'}}^U
  \end{gather*}

  \vspace{0.1in}
  Since  the above indices appearing  on both sides of the above equations are disjoint, both intersections on each side are transversal. This implies, in particular, that \\
  \begin{gather*}
  cod( \Delta_{J_1\cup J_2}^U,U)=cod(\Delta_{J_1}^U,U)+ cod(\Delta_{J_2}^U,U)  
  \end{gather*}

\vspace{0.1in}
which, in turn, implies that $|J_1\cup J_2|-1=(|J_1|-1)+(|J_2|-1)$, since, as we saw above, the diagonals are pure dimensional and the restrictions to $U$ are the unique components contained in $U$. This is a contradiction, since $J_1$ is disjoint from $J_2$.\\

Therefore, if $S\subseteq J_j$ for some $j\in\{1,2,\dots k'\}$, we have $\Delta_S^U \supseteq \Delta_{J_j}^U$, so we finally see that the minimal elements of $\wbuildrelU$ that contain $S_k=\Delta_{I_1}^U\cap \Delta_{I_2}^U\cap\dots \cap \Delta_{I_k}^U=\Delta_{J_1}^U\cap \Delta_{J_2}^U\cap\dots \cap \Delta_{J_{k'}}^U$ are precisely the elements  $\Delta_{J_j}^U$, where $j=1,2,\dots,k'$. This completes the proof of the claim.\\

Now, in view of the claim and since we wrote the elements of $\wbuildrel$ in ascending dimension order, we conclude, by \cite[Proposition 2.13]{Li}, that the sequence of blowups along the strict transforms of the subvarieties $\Delta_I^U\in \wbuildrelU$ coincides with the closure of the diagonal embedding \\
\begin{gather*}
U\setminus \bigcup\limits_ {\Delta_I^U \in \wbuildrelU} \Delta_I^U\xhookrightarrow\ U \times\displaystyle \prod \limits_ {\Delta_I^U\in\ \wbuildrelU} Bl_{\Delta_I^U} U
\end{gather*}

\vspace{0.1in}
hence the proposition.

   \qed
     
     \bigskip
\textbf{Remark}: In the above proof we assumed that the general fiber was reducible. In  particular, when $Y$ is an algebraically closed field, the above proposition shows that we can extend the definitions of section 2 to reducible schemes. 

\bigskip
\begin{thm}\label{thm7}
\begin{enumerate}
\item{The variety $\wmodrel$ is nonsingular. The boundary  $\wmodrel \setminus (X_Y^n \setminus \bigcup \limits_ {\Delta_I \in \wbuildrel} \Delta_I) $ is the union of $|\wbuildrel|$ divisors $D_Y^I,$ where $I\subset N   , |I|\geq2$ and $ \sum \limits _{i\in I} a_i >1\ $.}
\item{ Any set of boundary divisors intersects transversally. An intersection of divisors $D_Y^{I_1}\cap D_Y^{I_2}\cap\dots D_Y^{I_k}$ is nonempty precisely when the sets are nested in the sense that any pair $\{I_i,I_j\}$ either has empty intersection or one set is contained in the other. }
\end{enumerate}
\end{thm}
\textbf{Proof}: It suffices to prove the statement on the (open) irreducible components of $X_Y^n$. Then we have, from the proof of proposition \ref{prop6}, that the restrictions of the diagonals of  $\wbuildrel$ form a building set, so the proof is an immediate consequence of Theorem \ref{thm1}.\\

\textbf{Remark}: The divisors $D_Y^I$ above are precisely  the dominant transforms of the $\Delta_I\in\wbuildrel$, under the identification of $\wmodrel$ as a sequence of blowups $\wmodrel\rightarrow X_Y^n$ verified in proposition 6.\\

Now let us examine the fibers of the family $\wmodrel\rightarrow Y$. We have the following:\\
\begin{prop} \label{prop7}Let $f:X\rightarrow Y$ be a smooth morphism between smooth varieties. For a given $y\in Y$, let $X_y$ be the fiber of $f$ over $y$. Then there exists an isomorphism \\
\begin{gather*}
(\wmodrel)_y\iso (X_y)_{\mathcal{A}}[n]
\end{gather*}
\end{prop}
\textbf{Proof}: In Proposition \ref{prop6} we showed that the morphism $\wmodrel$ may be obtained from $X_Y^N$ as a sequence of blowups at centers that belong to $\wbuildrel$. The statement is Zariski local over the base, so after base changing to an irreducible component of $Y$, which is open since $Y$ is nonsingular, we may assume that the base is irreducible. We prove the proposition under the additional assumption that the general fiber is irreducible, which is the case of interest in the sequel, although it will be obvious from the proof that the general case is also true. Under this assumption, since the morphism $X_Y^N\rightarrow Y$ is smooth, so open, and because $Y$ is irreducible, we deduce that $X_Y^N$ is irreducible. The proof is  essentially identical to the proof of Proposition 1 in \cite{Pa1}. Let $m$ be the number of blowups needed and $\pi_m:\wmodrel \rightarrow X_Y^N$ be the resulting morphism. Also, denote by $X_Y^{N,i}$ be the variety resulting form the $i$-th blowup and let $\pi_i:X_Y^{N,i}\rightarrow X_Y^N$ the corresponding blowup morphism, where $0\leq i \leq m$. Finally let $\Delta_Y^i$ be the center of the ($i+1$)-th blowup $X_Y^{N,i+1}\rightarrow X_Y^{N,i}$ and let $\Delta_y^i$ be the center of the ($i$+1)-th blowup $X_y^{N,i+1}\rightarrow X_y^{N,i}$. Then we can show, by induction on the integer $i$, that there exists a natural isomorphism\\
\begin{gather}
(X_Y^{N,i})_y\iso X_y^{N,i}
\end{gather}
for all $y\in Y$. Here, $X_y^{N,i}$ is the variety obtained after blowing up the $i$-th element in the collection of subsets $\wbuild$ of $X_y^N$ (see section \ref{sec2.2}).

We need the following claim found in the proof of proposition 1 in \cite{Pa1}:\\

\textit{Claim}: Let $S$ be an irreducible smooth variety. Suppose that $R\rightarrow S$ is a smooth morphism and $T\xhookrightarrow\ R$ is a closed immersion such that $T$ is smooth over $S$. Then, for any $s\in S$ we have 
\begin{gather*}
(Bl_T R)_s\iso Bl_{T_s}R_s
\end{gather*}

\bigskip
When $i=0$, equation (16) reads $(X_Y^N)_y\iso X_y^N$, which is clear. Now, assume equation (16) holds for some $i$. From this we deduce the equality\\
\begin{gather}
(\Delta_Y^i)_y\iso \Delta_y^i
\end{gather}
where the above isomorphism is the restriction of the isomorphism in (16) to $(\Delta_Y^i)_y$. 
 $X_y^{N,i}$ is obtained from the smooth variety $X_y$ by sequence of blowups along smooth centers, so it is smooth. Consequently, we deduce from (16) that $(X_Y^{N,i})_y$ is also smooth for all $y\in Y$. This implies that the morphism\\
\begin{gather}
 X_Y^{N,i}\rightarrow Y 
\end{gather}
 is smooth. Moreover, we have seen that $\Delta_y^i$ is smooth by construction (section 2.2), so by (17), $(\Delta_Y^i)_y$ is also smooth, thus we conclude that the morphism\\
 \begin{gather}
 \Delta_Y^i\rightarrow Y
 \end{gather}
 is smooth. So, finally, we may apply the claim to  $\Delta_Y^i\xhookrightarrow\  X_Y^{N,i}$ and make use of the isomorphism (17) to conclude that the fiber of the blowup of  $X_Y^{N,i}$ at $\Delta_Y^i$ is isomorphic to the blowup of $X_y^N$ at $\Delta_y^i$, that is\\
 \begin{gather*}
 (X_Y^{N,i+1})_y\iso X_y^{N,i+1}
\end{gather*}
which proves the inductive statement (16) for $i+1$. By setting $i=m$ in (16) we complete the proof of proposition.\qed

\bigskip
\begin{subsubsection}{The relative weighted universal family} The construction of the universal family of Theorem \ref{thm3} can be generalized to the relative case  for our purposes. We construct the family as a sequence of blowups along subvarieties of $\wmodrel\times_Y X$.\\

 More specifically, consider the boundary divisors $D_{Y,I}$ of $\wmodrel$ (Theorem \ref{thm7}). Those are the iterated dominant transforms of the diagonals $\Delta_I$ of $\wbuildrel$; let $D_{Y,I}\rightarrow \Delta_I$ be the natural morphism; we may identify $\Delta_I$ with $X$ in $X^I$. Now consider the graph of the above morphism $D_{Y,I}\xhookrightarrow\ \wmodrel \times _Y X$ and define:
\begin{gather*}
D_{Y,I+}\subset \wmodrel\times_Y X
\end{gather*}
 to be its image. Now define the set:
 \begin{gather*}
 \wbuildrelu:=\{D_{Y,I+}\subset \wmodrel\times_Y X \, |\Delta_I \in \wbuildrel\}
  \end{gather*}
 and write its elements in order preserving bijection with the elements of $\Delta_I\in \wbuildrel$, i.e.  $D_{Y,I+}<D_{Y,J+}$ if and only if $\Delta_I<\Delta_J$. Note: we have seen that as long as the $\Delta_I$ are written in ascending dimension order, the resulting space $\wmodrel$ remains the same. We will see below that as long as this condition is satisfied, all  $D_{Y,I+}\in \wbuildrelu$ with $|I|=n-k$ become disjoint at the $k$-th step (see below). \\
 
We now define some intermediate steps $Z_k$ inductively as in section \ref{sec2.3}: To begin with, set
\begin{itemize}
\item $Z_0:=\wmodrel\times_Y X$. 
\item  $D_{Y,I+}^{(0)}:= D_{Y,I+}$, for all $I\subset N$ such that $\Delta_I\in  \wbuildrel$
\item  For each $i\in N$ we define $\Delta_{Y,i+}^{(0)}$ to be the embedding of $\wmodrel$ in $\wmodrel\times_Y X$ given by the graph of the morphism  $\wmodrel\rightarrow X_Y^n\xrightarrow{p_i} X$, where $p_i$ is the projection to the $i$-th factor.
\end{itemize}

 Assume $Z_k$ is defined. We will call $Z_k$ the \textbf{ $k$-th step}. Then we define: 
 \begin{itemize}
 \item$Z_{k+1}$ to be the blowup of $Z_k$ along the disjoint union of the subvarieties $D_{Y,I+}^{(k)}$, where $|I|=n-k$ and $\Delta_I \in \wbuildrel$; those subvarieties will be shown to be disjoint in the sequel. 
 \item for any $I\subset N$ s.t $\Delta_I\in\wbuild$, the variety $D_{Y,I+}^{(k+1)}$ is the dominant transform of $D_{Y,I+}^{(k)}$.
 \item for $i=1,\dots , n$, the variety $\Delta_{Y,i+}^{(k+1)}$ is the strict transform of $\Delta_{Y,i+}^{(k)}$.
   \end{itemize}
   
 Let $s$ be the greatest integer for which there exist a set $I\subset N$ such that $|I|=n-s$ and $\Delta_I\in \wbuildrel$. Then, the \textbf{universal family} is defined as\\
 \begin{center}
$\wmodrel^+:=Z_{s+1}$
\end{center}

\vspace{0.1in} 
and let 
\begin{center}
$\sigma_i:=\Delta_{Y,i+}^{(s+1)}$
\end{center}
for all $i=1,2,\dots,n$

As in Proposition \ref{prop7}, we can show that, if $y\in Y$ and $X_y$ is the fiber of $X\rightarrow Y$, then:\\
 \begin{gather}
  (\wmodrel^+)_y\iso (X_y)_{\mathcal{A}}[n]^+
  \end{gather}

 \bigskip
 Indeed, this is shown by induction on the number of steps $k$ as in Proposition 7. Let $|I|=n-k-1$. Also, suppose $V_k$ is the $k$-th step and $D_{I+}^{(k)}$ are the centers of the $k+1$-th step in the construction of  $(X_y)_{\mathcal{A}}[n]^+$, as defined in section \ref{sec2.3}. Finally let  ${\Delta}_{i+}^{(k)}$, $i=1,\dots n$, be the sections of $V_k\rightarrow (X_y)_{\mathcal{A}}[n]$. We can show that:\\
 \begin{gather*}
 {(Z_k)}_y\iso V_k
 \end{gather*}
 from which we deduce, via restriction, that
 \begin{gather}
    {(D_{Y,I+}^{(k+1)})}_y\iso D_{I+}^{(k+1)}
    \end{gather}
 and 
  \begin{gather}
    {(\Delta_{Y,i+}^{(k+1)})}_y\iso {\Delta}_{i+}^{(k+1)}
    \end{gather}

 When $k=0$, we have \\
 \begin{gather*}
 {(Z_0)}_y=(\wmodrel\times_Y X)_y=(\wmodrel_y)\times X_y\iso (X_y)_{\mathcal{A}}[n] \times X_y=V_0
    \end{gather*}
    by Proposition \ref{prop7} .   If  the claim is true for $k$, then, since all the centers in the construction of $(X_y)_{\mathcal{A}}[n]^+$ are smooth, we may apply the same argument as in Proposition \ref{prop7} to conclude that  \\
    \begin{gather*}
 {(Z_{k+1})}_y\iso V_{k+1}
 \end{gather*}

  hence the claim for $k+1$. This completes the induction.\\
  
   In particular, we see that the centers $D_{Y,I+}^{(k)}$ of the $k+1$-th step are pairwise disjoint for all $I\subset N$ such that $|I|=n-k-1$: indeed, we have already seen in proposition \ref{prop4}(ii) that the corresponding centers of the $k+1$-th step in the construction of $(X_y)_{\mathcal{A}}[n]^+$ are disjoint, so the result follows from equation (21).\\
   
   Now, let $\pi_\mathcal{A}:X_{Y}^{\mathcal{A}}[n]^+\rightarrow X_{Y}^{\mathcal{A}}[n]$ be the full resolution. Also, consider $y\in X_{Y}^{\mathcal{A}}[n]$ and let $y'$ be the image of $y$ in $Y$ under  $X_{Y}^{\mathcal{A}}[n]\rightarrow Y$. By proposition 7 and equations (20) and (22) we have a cartesian diagram over $y'\rightarrow Y$:\\
      \[
   \xymatrix
 {X_{y'}^{\mathcal{A}}[n] ^+\ar[d] \ar[r] &  X_{Y}^{\mathcal{A}}[n]^+\ar[d]\\
 X_{y'}^{\mathcal{A}}[n]\ar[r] \ar@/^/[u]^{\sigma'_i} & X_{Y}^{\mathcal{A}}[n]\ar@/^/[u]^{\sigma_i }}\\
  \]
  
  \bigskip
  From the above diagram we see that the fiber of $\pi_\mathcal{A}:X_{Y}^{\mathcal{A}}[n]^+\rightarrow X_{Y}^{\mathcal{A}}[n]$ over $y$ is the fiber of $X_{y'}^{\mathcal{A}}[n] ^+\rightarrow  X_{y'}^{\mathcal{A}}[n]$ over $y$, in other words a $n$-pointed $\mathcal{A}$ stable degeneration of $X_{y'}$. \\
  
  By the above discussion we conclude the following:
  \begin{prop}\label{prop8}The morphism $\pi_\mathcal{A}:X_{Y}^{\mathcal{A}}[n]^+\rightarrow X_{Y}^{\mathcal{A}}[n]$ along with the $\sigma_i,\,i=1,\dots,n$ is a universal family of $n$-pointed $\mathcal{A}$ stable degenerations on the fibers of $f:X\rightarrow Y$.\\
  \end{prop}
  \end{subsubsection}
  \end{subsection}

     \begin{subsection}{Relative weighted stable degenerations} \label{sec3.2}    Let $W$ be a smooth variety and $B$ a smooth curve with a distinguished point $b_0 \in B$. Also, let $\pi: W\rightarrow B$ be a flat morphism  such that the fibers $W_b$ are smooth for all $b \in B$ except  $b=b_0$, in which case $W_{b_0}$ is the union of two smooth varieties $X_1$ and $X_2$ intersecting transversally along a smooth divisor $D$. This implies that $N_{D/ X_1}\iso N_{D/X_2}^{\vee}$.    Now, let $\mathcal{A}:=\{a_1,a_2, \dots,a_n\} $ be a set of rational numbers such that $0<a_i\leq 1$ , for $i=1,2,\dots,n$. In this section we construct a compactification $W_{\pi,\mathcal{A}}[n] $ of the configuration space of $n$ labeled points in $W\rightarrow B$ equipped with weights $a_i$ such that 
     \begin{itemize}
 \item    the points lie in the smooth locus of $W\rightarrow B$
 \item for any set of labels $S\subset N$ of coincident points, we have $\sum \limits _{i\in S} a_i \leq 1\ $.
 \end{itemize}
 
First, we recall some preliminary definitions and results from \cite{AF2}. 
     \begin{mydef}(\cite[Definition 2.3.1] {ACFW})  The expansion of length $l\geq 0$ of $W_{b_0}$ is  the variety
     \begin{gather*}
     W(l):=X_1\bigsqcup\limits_{D=D_1^-} P_1\bigsqcup_{D_1^+=D_2^-}\cdots \bigsqcup_{D_{l-1}^+=D_l^-}P_l  \bigsqcup_{D_{l+}^+=D}X_2
               \end{gather*}
        where the $P_i$ are isomorphic to the $\PP^1$- bundle $\PP(N_{D/X_1}\oplus \mathcal{O}_D)=\PP( \mathcal{O}_D \oplus N_{D/X_2})$ and $D_j^-, D_j^+$ are the zero and infinity sections of $P_j$ respectively. The $P_j$ are called the exceptional components of the expansion.
     \end{mydef}
     \bigskip
    
     In \cite{AF2}, the authors construct a compactification of the space of $n$ tuples of points on the fibers of $\pi$ such that the points are not allowed to land in the singular locus, by replacing the critical fibers by expansions. We give the basic definitions and refer the reader to [ibid.] for more details. \\
     \begin{mydef}\cite[Definition 1.2.1]{AF2} A stable expanded configuration $(\mathcal{W},\widetilde{\sigma_i})$ of degree $n$ on $\pi:W\rightarrow B$ consists of:
     \begin{enumerate}
     \item either a smooth fiber $W_b$ or an expansion of $W_{b_0}$, which we denote by $\mathcal{W}$
     \item an ordered collection of $n$ smooth points $\widetilde{\sigma_i}\in \mathcal{W}^{sm},\,i=1,\dots n$
     \end{enumerate}
     such that in case $\mathcal{W}$ is an expansion, the following stability condition holds:
     \begin{itemize}
     \item each exceptional component $P_j$ contains at least one of the $\widetilde{\sigma_i}$.
     \end{itemize}
     \end{mydef}
     \bigskip
     There exists a \textit{moduli space of stable expanded configurations of degree $n$}, denoted by $W_{\pi}^{[n]}$, which is a smooth scheme by \cite[Proposition 1.2.5.]{AF2}. Also, there exisits a universal family $W_{\pi}^{[n]+}\rightarrow W_{\pi}^{[n]}$ together with $n$ sections $s_i:W_{\pi}^{[n]}\rightarrow W_{\pi}^{[n]+},\, i=1,\dots, n.$ We now give the following:
     
      \begin{mydef} An $n$-pointed relative $\mathcal{A}$ stable degeneration of $W\rightarrow B$ is pair $(\mathcal{W'},\sigma_i')$ such that 
      \begin{center}
    $W'=\mathcal{W}\bigsqcup\limits_{\mathcal{W}^{sm}\setminus\{\sigma_i\}} \widetilde{ \mathcal{W}}$
    \end{center}
    where $(\mathcal{W},\sigma_i)$ is a stable expanded configuration of degree $n$ and $(\widetilde{ \mathcal{W}},\sigma_i')$ is an $n$-pointed $\mathcal{A}$ stable degeneration of $(\mathcal{W}^{sm},\sigma_i)$ in the sense of section $2.2$.    
                  \end{mydef}    
                  
       We now construct the compactification of the configuration space described in the beginning of this section whose degenerate objects are the $n$-pointed relative $\mathcal{A}$ stable degenerations of $W\rightarrow B$. The resulting space is obtained from the moduli space of expanded configurations of degree $n$ \cite{AF2}, after pulling apart the sections as dictated by the weight data $\mathcal{A}$. \\
         
Let $W_{\pi}^{[n]}$ be the moduli space of stable expanded configurations. Since the sections $s_i:W_{\pi}^{[n]}\rightarrow W_{\pi}^{[n]+},\, i=1,\dots, n$  land in the relative smooth locus, $W_{\pi}^{[n]+,o}$, of $W_{\pi}^{[n]+}\rightarrow W_{\pi}^{[n]}$, we get $n$ morphisms $W_{\pi}^{[n]}\rightarrow W_{\pi}^{[n]+,o}$. For ease of notation set $X:=W_{\pi}^{[n]+,o}$ and $Y:=W_{\pi}^{[n]}$. Then, we have a morphism $Y\rightarrow X_Y^n$ induced by the $n$ sections and a natural morphism $\wmodrel\rightarrow X_Y^n$, where $\wmodrel$ is the relative weighted compactification defined in section 3.1. Define $W_{\pi,\mathcal{A}}[n]$ to be the fiber product :\\
  \[
\xymatrix{
W_{\pi,\mathcal{A}}[n] \ar[d] \ar[r] &X_{Y}^{\mathcal{A}}[n]\ar[d]\\
Y \ar[r]    	&X_Y^n }
\]
and let $W_{\pi,\mathcal{A}}[n]^+$ be the fiber product
\bigskip
  \[
   \xymatrix
 {W_{\pi,\mathcal{A}}[n]^+\ar[d] \ar[r] &  X_{Y}^{\mathcal{A}}[n]^+\ar[d]\\
W_{\pi,\mathcal{A}}[n] \ar[r] \ar@/^/[u]^{\sigma'_i} & X_{Y}^{\mathcal{A}}[n]\ar@/^/[u]^{\sigma_i }}\\
  \]
  
 where the sections $\sigma_i'$ are pullbacks of the sections $\sigma_i$. \\
 
By Propositions \ref{prop7} and \ref{prop8}, we see that $\wmodrel^+\rightarrow \wmodrel$ is a `universal family' of the $n$-pointed $\mathcal{A}$ stable degenerations on the fibers of $W_{\pi}^{[n]+,o}\rightarrow W_{\pi}^{[n]}$, 
so we conclude:
\begin{prop}\label{prop9}The morphism $W_{\pi,\mathcal{A}}[n]^+\rightarrow W_{\pi,\mathcal{A}}[n]$ along with the sections $\sigma_i'$ is a `universal family' of $n$-pointed relative $\mathcal{A}$ stable degenerations of $W\rightarrow B$.
\end{prop}

\end{subsection}
\end{section}
\begin{section}{Examples of weighted compactifications}\label{sec4}
\begin{subsection}{Weighted stable maps of degree one to $\PP^1$}In \cite{MM} the authors construct intermediate moduli spaces for the space of genus 0 stable maps to projective space. In the  case when the degree is 1 and the target is $\PP^1$, the following definition is given:\\
\begin{mydef}\cite[Definition 1.7] {MM}Let $k,n$ be integers such that $0\leq k\leq n$. Consider a morphism $\phi: C\rightarrow S\times \PP^1$ of degree 1 over each geometric fiber $C_s$ with $s\in S$ and $n$ marked sections of $C\rightarrow S$. The morphism will be called $k$- stable if the following conditions are satisfied for all $s\in S$:
\begin{itemize}
\item{not more than $n-k$ of the marked points in $C_s$ coincide};
\item{any end irreducible component of $C_s$, i.e. one that corresponds to a vertex of valence 1 in the dual graph of the curve, besides the parametrized one contains more than $n-k$ marked points };
\item{all marked sections lie in the relative smooth locus } ;
\item{$C_s$ has no non-trivial automorphisms fixing the points and the map to $\PP^1$}.
\end{itemize}
\end{mydef}  
Such morphisms are also called \textit{$k$ stable parametrized rational pointed curves}. The above moduli problem is finely represented by a scheme:

\bigskip
\begin{prop}\label{prop10}(\cite[Proposition 1.8]{MM}) 
For each $k$ such that $0\leq k\leq n$, let $\PP^1[n,k]$ be the closure of the diagonal embedding \\
\begin{align*}
(\PP^1)^n\setminus \bigcup_{i,j}\Delta_{i,j} \rightarrow (\PP^1)^n \times \prod_{|S|>n-k} Bl_{\Delta_S} (\PP^1)^S
\end{align*}

where $S\subset N=\{1,2,\dots,n\}$ and $|S|>2$. The moduli problem of $k$ stable parametrized rational pointed curves is finely represented by $\PP^1[n,k]$ . 
\end{prop}
\bigskip
The moduli space $\PP^1[n,k]$ is a weighted compactification. Indeed, take any rational number $\epsilon$ such that $\frac{1}{n-k+1}<\epsilon\leq \frac{1}{n-k}$ and consider the ordered set $\mathcal{A}_\epsilon:=\{\underbrace{\epsilon,\epsilon,\dots ,\epsilon}_{n\, elements}\}$. Then, the corresponding building set (section 2.2) is 
\begin{center}
$\mathcal{G}_{\mathcal{A}_\epsilon}=\{ \Delta_{12\dots n}; \Delta_{12\dots (n-1)},\dots ,\Delta_{23\dots n}; \dots;  \Delta_{12\dots n-k+1},\dots   ,\Delta_{k(k+1)\dots n} \}$, 
\end{center}
i.e. $\mathcal{G}_{\mathcal{A}_\epsilon}$ is the set of all diagonals of size greater than $n-k$ listed in ascending dimension order. By Definition 5, the weighted compactification, ${X_\mathcal{A_\epsilon}[n]}$, of $(\PP^1)^n\setminus \bigcup\limits_ {\Delta_S \in \mathcal{G}_{\mathcal{A_\epsilon}}} \Delta_S$ is the closure of the diagonal embedding\\

  \begin{center}
 $i: (\PP^1)^n\setminus\displaystyle \bigcup \limits_ {\Delta_S \in \mathcal{G}_{\mathcal{A_\epsilon}}} \Delta_S\xhookrightarrow\ (\PP^1)^n \times\displaystyle \prod \limits_ {\Delta_S\in \mathcal{G}_{\mathcal{A_\epsilon}}} Bl_{\Delta_S} (\PP^1)^n$
\end {center}

\vspace{0.1in}
Since the formation of blowups commutes with smooth base change, we have $Bl_{\Delta_S} (\PP^1)^n\iso {(\PP^1)^n} \times_{(\PP^1)^S}Bl_{\Delta_S} (\PP^1)^S  \iso (\PP^1)^{N\setminus S} \times Bl_{\Delta_S} (\PP^1)^S$. We now see that $(\PP^1)^n \times \displaystyle \prod\limits_{{\Delta_S \in \mathcal{G}_{\mathcal{A_\epsilon}}}} Bl_{\Delta_S} (\PP^1)^S$  may be realized as a closed subscheme of\\
\begin{center}
  $(\PP^1)^n \times\displaystyle \prod \limits_ {\Delta_S\in \mathcal{G}_{\mathcal{A_\epsilon}}} Bl_{\Delta_S} (\PP^1)^n=(\PP^1)^n  \times\displaystyle \prod \limits_ {\Delta_S\in \mathcal{G}_{\mathcal{A_\epsilon}}} (\PP^1)^{N\setminus S} \times\displaystyle \prod \limits_ {\Delta_S\in \mathcal{G}_{\mathcal{A_\epsilon}}} Bl_{\Delta_S} (\PP^1)^S$,
  \end{center} 
  via the graph of the morphism $(\PP^1)^n\rightarrow (\PP^1)^n \times\displaystyle \prod \limits_ {\Delta_S\in \mathcal{G}_{\mathcal{A_\epsilon}}} (\PP^1)^{N\setminus S}$, through which $i$ factors . Therefore the closure of $i$ is the same as the closure of the embedding $(\PP^1)^n\setminus\displaystyle \bigcup \limits_ {\Delta_S \in \mathcal{G}_{\mathcal{A_\epsilon}}} \Delta_S\xhookrightarrow\ (\PP^1)^n \times\displaystyle \prod \limits_ {\Delta_S\in \mathcal{G}_{\mathcal{A_\epsilon}}} Bl_{\Delta_S} (\PP^1)^S$. Now, since $\Delta_S\in \mathcal{G}_{\mathcal{A_\epsilon}}$ if and only if $|S|>n-k$ and $(\PP^1)^n\setminus \bigcup_{i,j}\Delta_{i,j}$ is an open subscheme of $(\PP^1)^n\setminus\displaystyle \bigcup \limits_ {\Delta_S \in \mathcal{G}_{\mathcal{A_\epsilon}}}\Delta_S\subset(\PP^1)^n$, we conclude that the closure of the embedding\\
\begin{align*}
(\PP^1)^n\setminus \bigcup_{i,j}\Delta_{i,j} \rightarrow (\PP^1)^n \times \prod_{|S|>n-k} Bl_{\Delta_S} (\PP^1)^S
\end{align*}
\bigskip
in Proposition 10 is the same as the closure of the image of $i$, that is  \\
\begin{center}
$\PP^1[n,k]\iso{X_\mathcal{A_\epsilon}[n]}$
\end{center}
\end{subsection}
\vspace{0.1in}
\end{section}
\begin{section}{The Chow ring of the weighted compactification}\label{sec5}Let  $X$ be a nonsingular variety over an algebraically closed field $k$. We express the Chow ring of the weighted compactification $\wmod$ as an algebra over the Chow ring of $X^n$. In \cite{R2} we use the results of this section to compute the Chow ring of Hassett's space in genus 0 \cite{Hassett}.
 In general, the Chow ring of the product $X^n$ can be hard to compute in terms of the Chow ring of $X$. For the remainder of this section we will assume that  $X$ has a \textit{cellular decomposition}, i.e. there exists a filtration\\
\begin{gather*}
X_n:=X\supset X_{n-1}\supset \dots X_0  \supset X_{-1}:=\emptyset 
\end{gather*}
such that each $X_i\setminus X_{i-1}$ is isomorphic to a disjoint union of affine spaces. The reason we make this assumption is that in this case there exists a natural isomorphism $A^{\bullet}(X^n)\iso (A^{\bullet}X)^{\otimes n}$. \\

To arrive at our result we need some preliminary notions and lemmas. Let $Y,Z$ be smooth varieties such that $Z$ is a codimension $d$ subvariety of $Y$, for some $d>0$. Assume, furthermore, that the restriction map of Chow rings $A^{\bullet}Y\rightarrow A^{\bullet}Z$ is surjective. Denote the kernel of the above map by $J_{Z/Y}$. With these assumptions, we have the following:
\begin{mydef}
A Chern polynomial for $Z\subset Y$, denoted by $P_{Z/Y}$, is a polynomial \\
\begin{gather*}
P_{Z/Y}=t^d+a_1t^{d-1}+\dots +a_{d-1}t +a_d \in (A^{\bullet}Y)[t]
\end{gather*}
where $a_i\in A^iY$ is a class restricting to $c_i(N_{Z/Y})$, the $i$-th Chern class of the normal bundle of $Z$ in $Y$. We further require that $a_d=[Z]$, the class of $Z$ in $ A^{\bullet}Y$.
\end{mydef}
\vspace{0.1in}
Note that Chern polynomials for $Z\subset Y$ are not unique. The fact that $A^{\bullet}Y\rightarrow A^{\bullet}Z$ is surjective guarantees the existence of a Chern polynomial for $Z\subset Y$.  \\

For our computations we will need a few lemmas:
\begin{lem}\label{lem3}\cite[Lemma 5.1 (b)]{FM}Let $V$ and $W$ be subvarieties of $Y$ meeting transversally in a variety $Z$. Suppose that $V$ and $W$ have Chern polynomials $P_{V/Y}(t)$ and $P_{W/Y}(t)$ respectively. Then $Z$ has a Chern polynomial
\begin{gather*}
P_{Z/Y}(t)=P_{V/Y}(t)\cdot P_{W/Y}(t)
\end{gather*}
\end{lem}
\end{section}
Let $\pi:\widetilde{Y}\rightarrow Y$ be the blowup of $Y$ along $Z$. We identify $A^{\bullet}Y$ as a subring of $A^{\bullet}\widetilde{Y}$ via $\pi^{*}:A^{\bullet}Y\rightarrow A^{\bullet}\widetilde{Y}$. Let $E$ denote the exceptional divisor. Then, we have the following:  
\begin{lem}\label{lem4}\cite[Lemma 5.2]{FM} Let $V$ be a smooth subvariety of $Y$ that is not contained in $Z$. Suppose there exists a Chern polynomial $P_{V/Y}(t)$ for $V\subset Y$. Then:
\begin{enumerate}
\item If $V$ meets $Z$ transversally, then $P_{V/Y}(t)$ is a Chern polynomial for $\widetilde{V}\subset\widetilde{Y}$.
\item If $V$ contains $Z$, then $P_{V/Y}(t-E)$ is a Chern polynomial for $\widetilde{V}\subset\widetilde{Y}$
\end{enumerate}
\end{lem}
Next, we have the following key formula due to Keel:
\begin{lem}\label{lem5}\cite{Keel} With the above assumptions and notation, the Chow ring of $\widetilde{Y}$ is given by\\
\begin{gather*}
A^{\bullet}\widetilde{Y}=\dfrac{A^{\bullet}Y[E]}{(J_{Z/Y}\cdot E, P_{Z/Y}(-E))}
\end{gather*}
\end{lem}
\bigskip
Finally, we record two more lemmas that allow us to compute Chow rings of blowups, again under the same hypotheses and notation as above.\\
\begin{lem}\label{lem6}\cite[Lemma 5.4]{FM} Suppose that $V$ is a nonsingular subvariety of $Y$ that intersects $Z$ transversally in a irreducible subvariety $V\cap Z$ and that the restriction  $A^{\bullet}V\rightarrow A^{\bullet}(V\cap Z)$ is surjective. Let $\widetilde{V}=Bl_Z V$. Then $A^{\bullet}\widetilde{Y}\rightarrow A^{\bullet}\widetilde{V}$ is surjective with kernel 
$J_{V/Y}$, if $V\cap Z$ is nonempty and kernel $(J_{V/Y}, E)$, if $V\cap Z$ is empty.
\end{lem}
\bigskip
\begin{lem}\label{lem7}\cite[Lemma 5.5]{FM} Suppose that $Z$ is the transversal intersection of nonsingular varieties $V$ and $W$ of $Y$ and that the restrictions $A^{\bullet}Y\rightarrow A^{\bullet}V$ and $A^{\bullet}Y\rightarrow A^{\bullet}W$ are also surjective. Then,
$A^{\bullet}\widetilde{Y}\rightarrow A^{\bullet}\widetilde{V}$ is surjective with kernel $(J_{V/Y}, P_{W/Y}(-E))$.
\end{lem}
\vspace{0.1in}
Now, let $S\subset N$ such that $\Delta_S\in \wbuild$. The restriction morphism $A^{\bullet}X^n\rightarrow A^{\bullet}\Delta_S$ is surjective. Define \\
\begin{gather*}
J_S:=ker(A^{\bullet}X^n\rightarrow A^{\bullet}\Delta_S)
\end{gather*}
By our assumptions, $X^n$ has a K\"{u}nneth decomposition, so we have\\
\begin{gather*}
J_S=\left(p_a^*x-p_b^*x\right)_{x\in A^{\bullet}X,\, a,b \in S}
\end{gather*}
where $p_a:X^n\rightarrow X$ is the projection onto the $a$-th factor. Also, define the polynomial $c_{a,b}(t) \in A^{\bullet}(X^n)[t]$ by \\
\begin{gather*}
c_{a,b}(t):=\displaystyle \sum\limits_{i=1}^{m}(-1)^i p_a^*(c_{m-i}) t^i +[\Delta_{\{a,b\}}]
\end{gather*}
where $m=dim (X)$ and $c_{m-i}$ is the $(m-i)$th Chern class of the tangent bundle $T_X$ of $X$. Finally, we have the following:\\
\begin{mydef}\begin{enumerate}
\item Let $S,T\subset N$ such that $|S|,|T|\geq2$. We say that $S$ and $T$ \textbf{overlap} if $S$ and $T$ have nonempty intersection and neither of them  is contained in the other. In this case we also say that $T$ is an \textbf{overlap} of $S$ and vice versa.
\item Let $S,T\subset N$ such that $|S|,|T|\geq2$. We say that $T$ is a \textbf{weak overlap} of $S$ (and vice versa) if the intersection of $S$ and $T$ is a singleton.
\end{enumerate}
\end{mydef}
\textbf{Remark}: In the above definition we do not require $\Delta_S$ or $\Delta_T$ to be in $\wbuild$.\\

We are now ready to state and prove the main result of this section: let $A^{\bullet}(X^n)[D^S]$ be the polynomial ring over $A^{\bullet}(X^n)$ in the variables $D^S$, where $S\subset N$ is such that $\Delta_S\in \wbuild$. This polynomial ring maps to the Chow ring of $\wmod$ by sending each $D^S$ to the class of the divisor $D_S$ of $\wmod$.\\
\begin{thm}\label{thm8}With the above assumptions and notation, the Chow ring of $\wmod$ is the quotient
\begin{gather*}
\displaystyle\dfrac{A^{\bullet}(X^n)[D^S]}{\mathcal{I}}
\end{gather*}
where $\mathcal{I}$ is the ideal generated by the following elements:
\begin{itemize}
\item $D^S\cdot D^T$, for any $S$ and $T$ that overlap and both $\Delta_S$ and $\Delta_T$ belong to $\wbuild$;\\
\item $J_S\cdot D^S$, where $\Delta_S\in \wbuild$;\\
\item $D^S\cdot \displaystyle \prod\limits_{k=1}^{|T|-1} c_{i_k,i_{k+1}}\left(\sum \limits_{I\supseteq S\cup T} D^I\right)$, where $T=\{i_1, i_2,\dots ,i_{|T|}\}$ is a weak overlap of $S$, $\Delta_S\in \wbuild$ and $\Delta_T$ is not necessarily in $\wbuild$;\\
\item$\displaystyle \prod\limits_{k=1}^{|S|-1} c_{i_k,i_{k+1}}\left(\sum \limits_{I\supseteq S} D^I\right)$, where $S=\{i_1,i_2,\dots ,i_{|S|}\}$ and $\Delta_S\in \wbuild$.
\end{itemize}
\end{thm}
\vspace{0.1in}
To prove the above theorem, we will successively apply Keel's formula (Lemma \ref{lem5}) to the blowups along the centers of $\wbuild$. In particular, we will need to know how the quantities $J_S$ and $P_{\Delta_I/X^n}$ change under the various blowups. We use the notation of Proposition \ref{prop1} to describe what happens in the intermediate steps:\\
 We extend the partial order of the set of indices of the elements of $\wbuild$ to a total order which respects inclusions of diagonals:
\begin{center}
 $I_0:=(12\dots n)<I_1<\dots <I_{|\wbuild|-1}$ \\
\end{center}
Define $X_0:=Bl_{{\Delta_{I_0}}}X^n$ and let $\wbuild^{(I_0)}:=\{ {\Delta_{I_0}}^{(I_0)} , {\Delta_{I_1}}^{(I_0)}, \dots {\Delta_{I_l}}^{(I_0)}, \dots ,{\Delta_{I_k}}^{(I_0)},\dots\}$ be the set of dominant transforms of the varieties(diagonals) that belong to $\wbuild$ in $X_0$. Then define $X_1:=Bl_{{\Delta_{I_1}}^{(I_{0})}}X_0$ and let  $\wbuild^{(I_1)}:=\{ {\Delta_{I_0}}^{(I_1)} , {\Delta_{I_1}}^{(I_1)}  \dots {\Delta_{I_l}}^{(I_1)}, \dots ,{\Delta_{I_k}}^{(I_1)},\dots\}$ be the set of dominant transforms of the varieties that belong to $\wbuild^{(I_0)}$ in $X_1$. For $k\geq2$, we may define inductively  the $I_k$-th iterated blowup, \\
\begin{center}
$X_k:= Bl_{\Delta_{I_{k}}^{(I_{k-1})}} X_{k-1}= Bl_{{\Delta_{I_k}}^{(I_{k-1})}}\dots Bl_{{\Delta_{I_0}}}X^n$
\end{center}
\vspace{0.1in}
and $\wbuild^{(I_k)}:=\{ {\Delta_{I_0}}^{(I_k)} ;\dots {\Delta_{I_l}}^{(I_k)}, \dots ,{\Delta_{I_k}}^{(I_k)},\dots\}$ as the set of dominant transforms of the subvarieties of $X_{k-1}$ that are the elements of $\wbuild^{(I_{k-1})}$. In other words, $\wbuild^{(I_k)}:=\{ {\Delta_{I_0}}^{(I_k)} ,\dots {\Delta_{I_l}}^{(I_k)}, \dots ,{\Delta_{I_k}}^{(I_k)},\dots\}$ is the set of the $I_k$-th iterated dominant transforms of the diagonals of $\wbuild$ in $X_k$. In general, for any $T\subset N$ such that $|T|\geq2$ and $\Delta_T$ \textbf{not necessarily in} $\wbuild$, we define ${\Delta_{T}}^{(I_k)}$ to be the iterated dominant transform of $\Delta_T$ under $X_k\rightarrow X^n$. \\

We have seen in the proof of Proposition \ref{prop1} that each $\wbuild^{(I_k)}$ is a building set. With the above notation we have the following: \\ 
\begin{lem}\label{lem8}Let $k>j\geq 0$ such that $I_j\supset I_k$. Consider an overlap $T$ of $I_k$ such that $I_k\cup T=I_j$. \\
\begin{enumerate}
\item If $l<j$ and $T$ is a weak overlap, then ${\Delta_{I_k}}^{(I_l)}$ and ${\Delta_T}^{(I_l)}$ intersect transversally and ${\Delta_{I_k}}^{(I_l)}\cap{\Delta_T}^{(I_l)}={\Delta_{I_j}}^{(I_l)}$.\\
\item ${\Delta_{I_k}}^{(I_j)}\cap{\Delta_T}^{(I_j)}=\emptyset$.
\end{enumerate}
\end{lem}

Now let $l\geq0$. Consider the set ${\wbuild}_{,l}$ of the \textbf{first $l+1$ elements} of $\wbuild$ with the above order. By Lemma \ref{lem1}, this is a building set. Now let $I\subset N$ such that $|I|\geq2$ and $I\subsetneq I_l$. We form the set $\mathcal{G}_{\mathcal{A},l}':={\wbuild}_{,l}\cup \{\Delta_I\}$ and equip its set of indices  with an order that extends the above total order of $\wbuild$ and respects inclusions of diagonals, i.e.
\begin{gather*}
I_0:=(12\dots n)<I_1<\dots <I_l<I
\end{gather*}
We can see that
\begin{lem}\label{lem9}
The set $\mathcal{G}_{\mathcal{A},l}':={\wbuild}_{,l}\cup \{\Delta_I\}$ is a building set.
\end{lem}
\textbf{Proof}: The proof is almost identical to the proof of Lemma \ref{lem1}, so we omit it.\qed 

\vspace{0.1in}
Also, as above, the set of iterated dominant transforms\\ 
\begin{center}
${\mathcal{G}_{\mathcal{A},l}'}^{(I_k)}:=\{ {\Delta_{I_0}}^{(I_k)} ,\dots {\Delta_{I_l}}^{(I_k)},{\Delta_{I}}^{(I_k)}\}$ 
\end{center}
\vspace{0.1in}
of the diagonals of $\mathcal{G}_{\mathcal{A},l}'$ 
in $X_k$ is a building set when $I_k\leq I$. \\

\textit{Proof of Lemma 8}:(1) We keep the above notation and prove the statement by induction. When $j=0$, the proof is immediate. Assume the statement is true for $j$. Let $l<j+1<k$ and consider a weak overlap $T$ of $I_k$ such that $I_k\cup T=I_{j+1}$. We have seen above that the set $\wbuild^{(I_{l-1})}$ of iterated dominant transforms of the diagonals in $X_{l-1}$  is a building set and that its elements either contain or intersect transversally with the center ${\Delta_{I_{l}}}^{(I_{l-1})}$. Using this, we distinguish the following cases:\\
\renewcommand{\theenumi}{\roman{enumi}}
\begin{enumerate}
\item
  If $I_{j+1}\subsetneq I_l$, then ${\Delta_{I_{j+1}}}^{(I_{l-1})}$ (and consequently each of the ${\Delta_{I_k}}^{(I_{l-1})}$ and ${\Delta_T}^{(I_{l-1})}$) contains ${\Delta_{I_l}}^{(I_{l-1})}$.
Since $l-1<j$, we have that ${\Delta_{I_k}}^{(I_{l-1})}$ and ${\Delta_T}^{(I_{l-1})}$ intersect transversally and ${\Delta_{I_k}}^{(I_{l-1})}\cap{\Delta_T}^{(I_{l-1})}={\Delta_{I_{j+1}}}^{(I_{l-1})}$. Then, by Lemma \ref{lem2}(3)(a), we conclude that ${\Delta_{I_k}}^{(I_{l})}$ and ${\Delta_T}^{(I_{l})}$ intersect transversally and ${\Delta_{I_k}}^{(I_{l})}\cap{\Delta_T}^{(I_{l})}={\Delta_{I_{j+1}}}^{(I_{l})}$, hence the statement for $j+1$.\\

\item If $I_{j+1}\not\subset I_l$, then  ${\Delta_{I_{j+1}}}^{(I_{l-1})}$ cannot contain the center so it has to be transversal to it.  Also, we have seen (see discussion after lemma 9) that the set  $\mathcal{G}_{\mathcal{A},{j+1}}':={\wbuild}_{,j+1}\cup \{\Delta_T\}$ is a building set (with order described above), where ${\wbuild}_{,j+1}$ is the set of the first $j+2$ elements of $\wbuild$. We have also seen that the set of the iterated dominant transforms ${\mathcal{G}_{\mathcal{A},j+1}'}^{(I_{l-1})}$ in $X_k$ of the diagonals that belong to $\mathcal{G}_{\mathcal{A},{j+1}}'$ is also a building set with the order described above. Then, since the order respects inclusions of diagonals, the center ${\Delta_{I_l}}^{(I_{l-1})}$ is minimal, so, by \cite[Lemma 2.6(i)]{Li}, we conclude that ${\Delta_T}^{(I_{l-1})}$ either contains the center or is transversal to it. Similarly, we can show that ${\Delta_{I_k}}^{(I_{l-1})}$ also has this property. All possible intersections of ${\Delta_{I_k}}^{(I_{l-1})}$ and ${\Delta_T}^{(I_{l-1})}$ with the center are then described as follows:\\
\begin{itemize}
\item Both ${\Delta_{I_k}}^{(I_{l-1})}$ and ${\Delta_T}^{(I_{l-1})}$ are transversal to the center ${\Delta_{I_l}}^{(I_{l-1})}$. Now, since $l-1<j$, we have that ${\Delta_{I_k}}^{(I_{l-1})}$ and ${\Delta_T}^{(I_{l-1})}$ intersect transversally and ${\Delta_{I_k}}^{(I_{l-1})}\cap{\Delta_T}^{(I_{l-1})}={\Delta_{I_{j+1}}}^{(I_{l-1})}$. Therefore, by Lemma \ref{lem2}(3)(b), we have that ${\Delta_{I_k}}^{(I_{l})}$ and ${\Delta_T}^{(I_{l})}$ intersect transversally and ${\Delta_{I_k}}^{(I_{l})}\cap{\Delta_T}^{(I_{l})}={\Delta_{I_{j+1}}}^{(I_{l})}$, hence the statement for $j+1$.\\

\item One of the ${\Delta_{I_k}}^{(I_{l-1})}$ and ${\Delta_T}^{(I_{l-1})}$ contains the center and the other is transversal to it. Let us assume, without loss of generality, that ${\Delta_{I_k}}^{(I_{l-1})}\supset{\Delta_{I_l}}^{(I_{l-1})}$ and ${\Delta_T}^{(I_{l-1})}$ is transversal to ${\Delta_{I_l}}^{(I_{l-1})}$. Since $l-1<j$, we have that ${\Delta_{I_k}}^{(I_{l-1})}$ and ${\Delta_T}^{(I_{l-1})}$ intersect transversally and ${\Delta_{I_k}}^{(I_{l-1})}\cap{\Delta_T}^{(I_{l-1})}={\Delta_{I_{j+1}}}^{(I_{l-1})}$. Then, by Lemma \ref{lem2}(3)(c), we have that ${\Delta_{I_k}}^{(I_{l})}$ and ${\Delta_T}^{(I_{l})}$ intersect transversally and ${\Delta_{I_k}}^{(I_{l})}\cap{\Delta_T}^{(I_{l})}={\Delta_{I_{j+1}}}^{(I_{l})}$, hence the statement for $j+1$.\\
\end{itemize}
\end{enumerate}
(2) Take a weak overlap $T'\subset T$ of $I_k$ such that $I_j=I_k\cup T'$. It suffices to show that   ${\Delta_{I_k}}^{(I_j)}\cap{\Delta_{T'}}^{(I_j)}=\emptyset$. By part (1), we have that ${\Delta_{I_k}}^{(I_{j-1})}$ and ${\Delta_{T'}}^{(I_{j-1})}$ intersect transversally and ${\Delta_{I_k}}^{(I_{j-1})}\cap{\Delta_{T'}}^{(I_{j-1})}={\Delta_{I_j}}^{(I_{j-1})}$. Using Lemma \ref{lem2}(2), we conclude the result. 
\qed\\

The main step in the proof of Theorem 8 is the following:
\begin{prop}\label{prop11}With notation as above:
\begin{enumerate}
\item  Let $I\subset N$ such that $|I|\geq2$ and either $\Delta_I\in \wbuild$ or $I\subsetneq I_j$; in the second case $\Delta_I$ need not necessarily be an element of $\wbuild$. Then\\
\begin{gather*}P_{\Delta_{I}^{(I_{j})}/X_j}(t)=P_{\Delta_I/X^n}(t-\sum\limits_{\substack{I_l\supseteq I\\I_l \leq I_j}} \Delta_{I_l}^{(I_j)})
\end{gather*} 
\item Let $j\geq 0$. For each $k$ such that $I_k>I_j$ (with the order of $\wbuild$ described above), the restriction $A^{\bullet}X_j\rightarrow A^{\bullet}{\Delta}_{I_k}^{(I_j)}$ is surjective with kernel J$_{I_k}^{\,\,(j)}$ generated by\\
\begin{itemize}
\item $J_{I_k}$;\\
\item $\Delta_{I_m}^{(I_j)}$, for all $I_m$ that overlap with $I_k$ such that $I_m\leq I_j$;\\
\item $P_{\Delta_{I}/X^n}(-\sum\limits_{\substack{I_l\supseteq I\\I_l \leq {I\cup I_k}}} \Delta_{I_l}^{(I_j)})$, for all weak overlaps $I$ of $I_k$ such that $I\cup I_k\leq I_j$.
\end{itemize}
\end{enumerate}
\end{prop} 
\textbf{Remark}: In the third bullet of (2), $\Delta_{I\cup I_k}\in \wbuild$, since $I_k \cup I$ contains $I_k$, so $I\cup I_k\leq I_j$ makes sense. \\

\textbf{Proof}:(1) Let $j\geq0$ and consider $\mathcal{G}_{\mathcal{A},j}':={\wbuild}_{,j}\cup \{\Delta_I\}$, where, as above, ${\wbuild}_{,j}$ is the set of the first $j+1$ elements of $\wbuild$. The assumption of the Proposition on $I$ is equivalent to saying that $I_j<I$ with respect to the order of $\mathcal{G}_{\mathcal{A},j}'$ described above. We proceed by induction on $j$. When $j=0$ the statement reads 
$P_{\Delta_{I}^{(I_{0})}/X_0}(t)=P_{\Delta_I/X^n}(t-\Delta_{I_0}^{(I_0)})$. By the hypothesis $I\subsetneq I_0:=(12\dots n)$, so $\Delta_I\supsetneq \Delta_ {I_0}$, hence the result follows from Lemma \ref{lem4}(2). \\

Now, we assume that the statement is true for $j$. We want to prove the statement for $j+1$. Let $\mathcal{G}_{\mathcal{A},j+1}'={\wbuild}_{,j+1}\cup \{\Delta_I\}$ where ${\wbuild}_{,j+1}$ is the set of the first $j+1$ elements of $\wbuild$. As noted above, the condition of the proposition for $j+1$ is equivalent to $I_{j+1}<I$ (therefore $I_j<I$). Now the statement for $j$ reads\\
\begin{gather}P_{\Delta_{I}^{(I_{j})}/X_j}(t)=P_{\Delta_I/X^n}(t-\sum\limits_{\substack{I_l\supseteq I\\I_l \leq I_j}} \Delta_{I_l}^{(I_j)})
\end{gather} 
 We saw above that the set 
\begin{center}
${\mathcal{G}_{\mathcal{A},j+1}'}^{(I_j)}:=\{ {\Delta_{I_0}}^{(I_{j})} ,\dots, {\Delta_{I_{j}}}^{(I_{j})}, {\Delta_{I_{j+1}}}^{(I_{j})},{\Delta_{I}}^{(I_{j})}\}$ 
\end{center}
of the iterated dominant transforms of the varieties that belong to $\mathcal{G}_{\mathcal{A},j+1}'$ is a building set and that the center ${\Delta_{I_j+1}}^{(I_{j})}$ is minimal (see proof of Proposition \ref{prop1}). Then, by \cite[Lemma 2.6(i)]{Li} the rest of the varieties in $\mathcal{G}_{\mathcal{A},j+1}'$ either contain or are transversal to ${\Delta_{I_{j+1}}}^{(I_j)}$. In particular, if $I\subsetneq I_{j+1}$  then ${\Delta_{I}}^{(I_j)}\supseteq {\Delta_{I_{j+1}}}^{(I_j)}$, otherwise ${\Delta_{I}}^{(I_j)}$ is transversal to ${\Delta_{I_{j+1}}}^{(I_j)}$. Therefore, by Lemma \ref{lem4} we have\\
\begin{gather}
P_{\Delta_{I}^{(I_{j+1})}/X_{j+1}}(t)=
\begin{cases}
P_{\Delta_{I}^{(I_{j})}/X_{j}}(t-\Delta_{I_{j+1}}^{(I_{j+1})}),\,\text{if}\, I\subsetneq I_{j+1}\\[1em]
P_{\Delta_{I}^{(I_{j})}/X_{j}}(t),\, \text{otherwise}
\end{cases}
\end{gather}

\vspace{0.1in}

Now, by (23) and (24) and since the divisors  $\Delta_{I_l}^{(I_j)}$ pull back to $\Delta_{I_l}^{(I_{j+1})}$ in $X_{j+1}$, we deduce that \\
\begin{gather}
P_{\Delta_{I}^{(I_{j+1})}/X_{j+1}}(t)=
\begin{cases}
P_{\Delta_I/X^n}(t-\Delta_{I_{j+1}}^{(I_{j+1})} -\sum\limits_{\substack{I_l\supseteq I\\I_l \leq I_j}} \Delta_{I_l}^{(I_{j+1})}),\,\text{if}\, I\subsetneq I_{j+1}\\[2em]
P_{\Delta_I/X^n}(t-\sum\limits_{\substack{I_l\supseteq I\\I_l \leq I_j}} \Delta_{I_l}^{(I_{j+1})}),\, \text{otherwise}
\end{cases}
\end{gather}
in other words \\
\begin{gather*}
 P_{\Delta_{I}^{(I_{j+1})}/X_{j+1}}(t)=P_{\Delta_I/X^n}(t-\sum\limits_{\substack{I_l\supseteq I\\I_l \leq I_{j+1}}} \Delta_{I_l}^{(I_{j+1})})
 \end{gather*}
 hence the statement for $j+1$ and this completes the induction.\\
 
 (2) We proceed again by induction: let $j=0$. For any $k>0$ (i.e $I_k>I_0$), we have $I_k\subsetneq I_0:=(12\dots n)$. Consider a weak overlap $I$ of $I_k$ such that $I_0=I_k\cup I$. The intersection $\Delta_{I_k}\cap \Delta _I=\Delta_{I_0}$ is transversal, so by Lemma \ref{lem7}, we have that $A^{\bullet}X_0\rightarrow A^{\bullet}{\Delta}_{I_k}^{(I_0)}$ is surjective with kernel $J_{{I_k}}^{(0)}=(J_{I_k}, P_{\Delta_I/X^n}(-\Delta_{{I_0}}^{(I_0)}))$, which is the statement of the proposition for $j=0$.\\
 
Now assume that the statement is true for $j$ and consider any $I_k>I_{j+1}$. We distinguish three cases:
\renewcommand{\theenumi}{\roman{enumi}}
\begin{enumerate}
\item $I_k\subset I_{j+1}$: in this case, $\Delta_{I_k}^{(I_j)}$ contains the center, $\Delta_{I_{j+1}}^{(I_j)}$, of the blowup $X_{j+1}$. Now, let $I$ be any weak overlap of $I_k$ such that $I_{j+1}=I_k\cup I$.  By Lemma \ref{lem8}(1), we have ${\Delta_{I_k}}^{(I_j)}\cap{\Delta_I}^{(I_j)}={\Delta_{I_{j+1}}}^{(I_j)}$, where the intersection ${\Delta_{I_k}}^{(I_j)}\cap{\Delta_I}^{(I_j)}$ is transversal. Then, by Lemma \ref{lem7}, we see that the restriction $A^{\bullet}X_{j+1}\rightarrow A^{\bullet}{\Delta}_{I_k}^{(I_{j+1})}$ is surjective with kernel $J_{I_k}^{\,\,(j+1)}$ generated by $J_{I_k}^{\,\,(j)}$ and $P_{\Delta_{I}^{(I_{j})}/X_j}(-\Delta_{I_{j+1}}^{(I_{j+1})})$. Since we assumed the claim is true for $j$ and since the divisors $\Delta_{I_l}^{(I_j)}$ and $\Delta_{I_m}^{(I_j)}$ pullback to $\Delta_{I_l}^{(I_{j+1})}$ and  $\Delta_{I_m}^{(I_{j+1})}$ respectively in $X_{j+1}$, the contribution from $J_{I_k}^{\,\,(j)}$ is:\\
\begin{itemize}
\item $J_{I_k}$;\\
\item $\Delta_{I_m}^{(I_{j+1})}$, for all $I_m$ that overlap with $I_k$ such that $I_m\leq I_j$;\\
\item $P_{\Delta_{I}/X^n}(-\sum\limits_{\substack{I_l\supseteq I\\I_l \leq {I\cup I_k}}} \Delta_{I_l}^{(I_{j+1})})$, for all weak overlaps $I$ of $I_k$ such that $I\cup I_k\leq I_j$.
\end{itemize}
\vspace{0.1in}
Also, by part (1) above, we see, since  the divisors $\Delta_{I_l}^{(I_j)}$ pull back to the divisors $\Delta_{I_l}^{(I_{j+1})}$ in $X_{j+1}$, that the contribution of $P_{\Delta_{I}^{(I_{j})}/X_j}(-\Delta_{I_{j+1}}^{(I_{j+1})})$ is:
\begin{gather*}
P_{\Delta_{I}^{(I_{j})}/X_j}(-\Delta_{I_{j+1}}^{(I_{j+1})})=P_{\Delta_I/X^n}(-\sum\limits_{\substack{I_l\supseteq I\\I_l \leq I_j}} \Delta_{I_l}^{(I_{j+1})}-\Delta_{I_{j+1}}^{(I_{j+1})})=P_{\Delta_I/X^n}(-\sum\limits_{\substack{I_l\supseteq I\\I_l \leq I_{j+1}}} \Delta_{I_l}^{(I_{j+1})}). \\
\end{gather*}
The combination of these contributions proves the statement for $j+1$ in case (i).\\

\item $I_{j+1}$ is an overlap of $I_k$: First, as we have seen above, ${\Delta_{I_k}}^{(I_j)}$ either contains or intersects transversally with the center of the blowup $X_{j+1}$, which is $\Delta_{I_{j+1}}^{(I_j)}$. Since $I_{j+1}$ is an overlap of $I_k$, we can eliminate the first case, so ${\Delta_{I_k}}^{(I_j)}$  must intersect transversally with $\Delta_{I_{j+1}}^{(I_j)}$. Also, we can show that ${\Delta_{I_k}}^{(I_j)}\cap \Delta_{I_{j+1}}^{(I_j)}=\emptyset$. Indeed, observe that $I_{j+1}\cup I_k \in \wbuild$ since it contains $I_{j+1}\in \wbuild$. Then, we can see that  $I_{j+1}\cup I_k\leq I_j$: indeed, since 
$I_{j+1}$ is an overlap of $I_k$, we see that $|I_{j+1}\cup I_k|>|I_{j+1}|$, which implies that  $I_{j+1}\cup I_k<I_{j+1}$, therefore $I_{j+1}\cup I_k\leq I_j$. In case $I_{j+1}\cup I_k= I_j$ , we conclude from Lemma \ref{lem8}(2) that ${\Delta_{I_k}}^{(I_j)}\cap \Delta_{I_{j+1}}^{(I_j)}=\emptyset$. In case $I_{j+1}\cup I_k< I_j$, set $I_{j'}:=I_{j+1}\cup I_k$ (note $j'<j$). We may apply Lemma \ref{lem8}(2) again to deduce that ${\Delta_{I_k}}^{(I_{j'})}$ and $\Delta_{I_{j+1}}^{(I_{j'})}$ are disjoint, which implies that the strict transforms will still be disjoint after the sequence of blowups $X_j\rightarrow X_{j'}$, that is ${\Delta_{I_k}}^{(I_j)}\cap \Delta_{I_{j+1}}^{(I_j)}=\emptyset$.\\

We may now apply Lemma 6 to deduce that the restriction $A^{\bullet}X_{j+1}\rightarrow A^{\bullet}{\Delta}_{I_k}^{(I_{j+1})}$ is surjective with kernel $J_{I_k}^{\,\,(j+1)}$ generated by 
\begin{enumerate}
\item $J_{I_k}^{\,\,(j)}$ and 
\item the exceptional divisor $\Delta_{I_{j+1}}^{(I_{j+1})}$. 
\end{enumerate}
Since we assumed the claim is true for $j$ and since the divisors $\Delta_{I_l}^{(I_j)}$ and $\Delta_{I_m}^{(I_j)}$ pullback to $\Delta_{I_l}^{(I_{j+1})}$ and  $\Delta_{I_m}^{(I_{j+1})}$ respectively in $X_{j+1}$, the contribution from $J_{I_k}^{\,\,(j)}$ is:\\
\begin{itemize}
\item $J_{I_k}$;\\
\item $\Delta_{I_m}^{(I_{j+1})}$, for all $I_m$ that overlap with $I_k$ such that $I_m\leq I_j$;\\
\item $P_{\Delta_{I}/X^n}(-\sum\limits_{\substack{I_l\supseteq I\\I_l \leq {I\cup I_k}}} \Delta_{I_l}^{(I_{j+1})})$, for all weak overlaps $I$ of $I_k$ such that $I\cup I_k\leq I_j$.
\end{itemize}
\vspace{0.1in}
Observe that there can be no weak overlap $I$ such that $I\cup I_k= I_{j+1}$, so the above contributions combined with the contribution $\Delta_{I_{j+1}}^{(I_{j+1})}$ are precisely those mentioned in the statement of the proposition for $j+1$.\\
\item $I_{j+1}$ is disjoint from $I_k$: In this case ${\Delta_{I_k}}^{(I_j)}\cap \Delta_{I_{j+1}}^{(I_j)}\neq\emptyset$ and the result follows from Lemma \ref{lem6} as above. 
\end{enumerate}
\qed

\bigskip
We now compute the Chow ring of each of the intermediate varieties $X_j$. Once again, we keep the above notation.\\
\begin{prop}\label{prop12} Let $j\geq0$. The Chow ring of $X_j$ is the quotient \\
\begin{gather*}
\displaystyle\dfrac{A^{\bullet}(X^n)[\{{\Delta_{I_i}}^{(I_j)}\}_{i\leq j}]}{\mathcal{I}_j}
\end{gather*}
where $\mathcal{I}_j$ is the ideal generated by the following elements:\\
\begin{itemize}
\item ${\Delta_{I_k}}^{(I_j)}\cdot {\Delta_{I_l}}^{(I_j)}$ for any $I_k, I_l\leq I_j$ that overlap;\\
\item $J_{I_k}\cdot {\Delta_{I_k}}^{(I_j)}$ for all $I_k\leq I_j$;\\
\item ${\Delta_{I_k}}^{(I_j)}\cdot P_{\Delta_{I}/X^n}(-\sum\limits_{\substack{I_l\supseteq I\\I_l \leq {I\cup I_k}}} \Delta_{I_l}^{(I_j)})$ for all $I_k\leq I_j$ and all weak overlaps $I$ of $I_k$ such that $I_k\cup I\leq I_{j-1}$;\vspace{0.1in}

\item $P_{\Delta_{I_k}/X^n}(-\sum\limits_{\substack{I_l\supseteq I_k\\I_l \leq I_j}} \Delta_{I_l}^{(I_j)})$ for all $I_k\leq I_j$.
\end{itemize}
\end{prop}
\textbf{Proof}: We prove the proposition by induction on $j$. Let $j=0$; then, by Lemma 5 the Chow ring of $X_0=Bl_{\Delta_{I_0}}X^n$ is given by \\
\begin{gather*}
A^{\bullet}X_0=\dfrac{(A^{\bullet}X^n)[{\Delta_{I_0}}^{(I_0)}]}{(J_{I_0}\cdot {\Delta_{I_0}}^{(I_0)}, P_{\Delta_{I_0}/X^n}(-{\Delta_{I_0}}^{(I_0)}))}
\end{gather*}
which proves the statement of the proposition for $j=0$.\\

Now assume that the statement is true for some $j\geq0$. Again, by Keel's formula (Lemma \ref{lem5}), we conclude that 
$A^{\bullet}X_{j+1}$ is the quotient of $(A^{\bullet}X_j)[{\Delta_{I_{j+1}}}^{(I_{j+1})}]$ by the ideal\\
\begin{gather*}
 \mathcal{I}_{j+1,j}:=(J_{I_{j+1}}^{\,(j)}\cdot {\Delta_{I_{j+1}}}^{(I_{j+1})}, P_{\Delta_{I_{j+1}}^{(I_j)}/X_j}(-{\Delta_{I_{j+1}}}^{(I_{j+1})}))
\end{gather*}
\bigskip
Applying Proposition \ref{prop11} for $k=j+1$, since the divisors $\Delta_{I_l}^{(I_j)}$ and $\Delta_{I_m}^{(I_j)}$ pullback to $\Delta_{I_l}^{(I_{j+1})}$ and  $\Delta_{I_m}^{(I_{j+1})}$ respectively in $X_{j+1}$, we deduce that the contribution from $J_{I_{j+1}}^{\,\,(j)}\cdot {\Delta_{I_{j+1}}}^{(I_{j+1})}$ is:\\
\begin{itemize}
\item $J_{I_{j+1}}\cdot {\Delta_{I_{j+1}}}^{(I_{j+1})}$;\\
\item $\Delta_{I_m}^{(I_{j+1})}\cdot {\Delta_{I_{j+1}}}^{(I_{j+1})}$, for all $I_m$ that overlap with $I_{j+1}$ such that $I_m\leq I_j$;\\
\item $ {\Delta_{I_{j+1}}}^{(I_{j+1})}\cdot P_{\Delta_{I}/X^n}(-\sum\limits_{\substack{I_l\supseteq I\\I_l \leq {I\cup I_{j+1}}}} \Delta_{I_l}^{(I_{j+1})})$, for all weak overlaps $I$ of $I_{j+1}$ such that $I\cup I_{j+1}\leq I_{j}$.
\end{itemize}
\vspace{0.1in}
Also, the contribution from $P_{\Delta_{I_{j+1}}^{(I_j)}/X_j}(-{\Delta_{I_{j+1}}}^{(I_{j+1})}))$ is \\
\begin{center}
$P_{\Delta_{I_j+1}/X^n}(-\sum\limits_{\substack{I_l\supseteq I_{j+1}\\I_l \leq I_j}} \Delta_{I_l}^{(I_{j+1})}-\Delta_{I_{j+1}}^{(I_{j+1})})=P_{\Delta_{{I_j+1}/X^n}}(-\sum\limits_{\substack{I_l\supseteq I_{j+1}\\I_l \leq I_{j+1}}} \Delta_{I_l}^{(I_{j+1})})$ 
\end{center}
due to Proposition 11(1) and the fact that the divisors $\Delta_{I_l}^{(I_j)}$ pull back to $\Delta_{I_l}^{(I_{j+1})}$.\\
 
All these contributions give the total contribution of $\mathcal{I}_{j+1,j}$ above. Therefore, we see immediately that the total contribution of $\mathcal{I}_{j+1}$, which is the contribution of $\mathcal{I}_{j+1,j}$ along with the contribution of $\mathcal{I}_{j}$, is precisely the set of relations described in the
statement of the proposition for $j+1$. This completes the induction.\qed\\

\vspace{0.1in}

Now the proof of Theorem 8 is an immediate consequence of the above Proposition for $j=|\wbuild|-1$ and the following:\\
\begin{lem}\label{lem10} For any $I\subset N$ such that $|I|\geq2$, the Chern polynomial $P_{\Delta_I/X^n}(t)$ is equal to the product of Chern polynomials 
$\displaystyle \prod\limits_{k=1}^{|I|-1} c_{i_k,i_{k+1}}(-t)$, where $I=\{i_1,i_2,\dots ,i_{|I|}\}$.
\end{lem}
\textbf{Proof}: We prove this by induction on $|I|$. When $|I|=2$, the result is immediate. Let $j\geq3$ and assume the result is true for all $I\subset N$ with $|I|<j$. Then consider a set with $j$ elements $I':=\{i_1,i_2,\dots i_j\}$. Set $J:=\{i_1,i_2\}$ and $J' :=\{i_2,\dots i_j\}$. Then the intersection $\Delta_{J}\cap\Delta_{J'}=\Delta_{I'}$ is transversal. By Lemma \ref{lem3}, we have\\
\begin{center}
 $P_{\Delta_{I'}/X^n}(t)=P_{\Delta_J/X^n}(t)\cdot P_{\Delta_{J'}/X^n}(t)$
\end{center}
\vspace{0.1in}
which, by the induction hypothesis is equal to\\ 
\begin{gather*}
c_{i_1,i_2}(-t)\cdot \displaystyle \prod\limits_{k=2}^{j-1} c_{i_k,i_{k+1}}(-t)=\displaystyle \prod\limits_{k=1}^{|I'|-1} c_{i_k,i_{k+1}}(-t)
\end{gather*}
hence the statement for $j+1$.\qed
\bigskip

 \bibliographystyle{amsalpha}             
\bibliography{auta}
\end{document}